\documentclass[11pt]{article}
\usepackage{color}
\usepackage{amssymb}
\usepackage{amsthm,array,amssymb,amscd,amsfonts,latexsym, url}
\usepackage{amsmath}
\newtheorem{theo}{Theorem}[section]
\newtheorem{sublemm}[theo]{Sublemma}

\newtheorem{prop}[theo]{Proposition}
\newtheorem{lemm}[theo]{Lemma}
\newtheorem{coro}[theo]{Corollary}
\newtheorem{rema}[theo]{Remark}
\newtheorem{Defi}[theo]{Definition}

\newtheorem{conj}[theo]{Conjecture}
\newtheorem{question}[theo]{Question}

\newcommand{\cqfd}
{%
\mbox{}%
\nolinebreak%
\hfill%
\rule{2mm}{2mm}%
\medbreak%
\par%
}
\newfont{\gothic}{eufb10}

\date{}

\begin{document}
\title{Chow rings and  decomposition theorems for  families of  $K3$ surfaces and Calabi-Yau hypersurfaces}
\author{Claire Voisin
\\CNRS, Institut de math\'{e}matiques de Jussieu}

 \maketitle \setcounter{section}{0}

\begin{abstract}  The decomposition theorem for smooth projective morphisms
$\pi:\mathcal{X}\rightarrow B$ says that
$R\pi_*\mathbb{Q}$ decomposes as $\oplus R^i\pi_*\mathbb{Q}[-i]$. We describe simple examples
  where it is not  possible to have such a decomposition
compatible with cup-product, even after restriction to Zariski dense open sets of $B$. We prove however that  this is always possible for families of $K3$ surfaces (after
shrinking the base), and
show how this result relates to a result by Beauville and the author (\cite{beauvoi}) on the Chow ring of $K3$ surfaces
$S$. We give two proofs of this result, the first one involving
$K$-autocorrespondences of $K3$ surfaces, seen as analogues of isogenies of abelian varieties, the second one involving a
certain decomposition of the small diagonal in $S^3$ obtained in \cite{beauvoi}.
We also prove an analogue of such a decomposition of the small diagonal in $X^3$ for Calabi-Yau hypersurfaces $X$
in $\mathbb{P}^n$, which in turn
provides strong restrictions on their Chow ring.
\end{abstract}

Let $\pi:\mathcal{X}\rightarrow B$ be  a smooth projective morphism. The decomposition theorem, proved by
Deligne in \cite{deligne} as a consequence of the hard Lefschetz theorem, is the following statement:
\begin{theo} (Deligne 1968) In the derived category of sheaves of $\mathbb{Q}$-vector spaces on $B$, there is a decomposition
\begin{eqnarray}
\label{isodec}
R\pi_*\mathbb{Q}=\oplus_i R^i\pi_*\mathbb{Q}[-i].
\end{eqnarray}
\end{theo}
This statement is  equivalent, as explained by Deligne in {\it loc. cit.} to (a universal version of)
the degeneracy at $E_2$ of the Leray spectral sequence of $\pi$. Deligne came back in   \cite{delignemotivedec} to the problem
of constructing a {\it canonical} such decomposition, given the topological Chern class $l$ of a
relatively ample line bundle on $\mathcal{X}$ and imposing partial compatibilities with the morphism of  cup-product
with $l$.

Note that both sides of (\ref{isodec}) carry a cup-product. On the right, we put the
direct sum of the cup-product maps $\mu_{i,j}:R^i\pi_*\mathbb{Q}\otimes R^j\pi_*\mathbb{Q}\rightarrow R^{i+j}\pi_*\mathbb{Q}$.
On the left, one needs to choose an explicit representation of $R\pi_*\mathbb{Q}$ by a complex
$C^*$, together with
an explicit morphism of complexes $\mu:C^*\otimes C^*\rightarrow C^*$ which induces the cup-product in cohomology.
 When passing to
coefficients $\mathbb{R}$ or $\mathbb{C}$, one can take $C^*=\pi_*\mathcal{A}^*_\mathcal{X}$, where
$\mathcal{A}^*_\mathcal{X}$ is the sheaf of $\mathcal{C}^\infty$ real or complex differential forms on $\mathcal{X}$
 and for $\mu$ the wedge product of forms. For rational coefficients, the explicit construction
 of the cup-product at the level of complexes (for example \v{C}ech complexes)
 is more painful (see \cite[6.3]{godement}).
 The resulting cup-product  morphism $\mu$ will be canonical only in the derived category.

The question we study in this paper is the following:
\begin{question}\label{question} Given a family of smooth projective varieties $\pi:\mathcal{X}\rightarrow B$,  does there exist a decomposition
as above  which is multiplicative, that is compatible with the morphism
$$\mu:R\pi_*\mathbb{Q}\otimes R\pi_*\mathbb{Q}\rightarrow R\pi_*\mathbb{Q}$$
given by cup-product?
\end{question}

Let us give three examples: In the first one, which is the case of families of abelian varieties, the answer to  Question
 \ref{question} is affirmative. This was proved by Deninger and Murre in \cite{deningermurre} as a consequence of a
 much more general ``motivic'' decomposition result.
\begin{prop} \label{abelian}  For any family $\pi:\mathcal{A}\rightarrow B$ of abelian varieties (or complex tori), there is
a multiplicative decomposition isomorphism $R\pi_*\mathbb{Q}=\oplus_i R^i\pi_*\mathbb{Q}[-i]$.
\end{prop}

In the next two examples, the answer to Question \ref{question} is negative.
The simplest  example is that of projective bundles $\pi:\mathbb{P}(\mathcal{E})\rightarrow B$, where
$\mathcal{E}$ is a locally free sheaf on $B$.
\begin{prop} \label{lemmac1} Assume that $c_1^{top}(\mathcal{E})=0$ in $H^2(B,\mathbb{Q})$. Then, if there exists a multiplicative decomposition
isomorphism for $\pi:\mathbb{P}(\mathcal{E})\rightarrow B$, one  has  $c_i^{top}(\mathcal{E})=0$ in $H^{2i}(B,\mathbb{Q})$ for all $i>0$.
\end{prop}
{\bf Proof.} Let $h=c_1^{top}(\mathcal{O}_{\mathbb{P}(\mathcal{E})}(1))\in H^2(\mathbb{P}(\mathcal{E}),\mathbb{Q})$.
It is standard that  $$H^2(\mathbb{P}(\mathcal{E}),\mathbb{Q})=\pi^*H^2(B,\mathbb{Q})\oplus \mathbb{Q}h,$$ where
$\pi^*H^2(B,\mathbb{Q})$ identifies canonically with the deepest term
$H^2(B,R^0\pi_*\mathbb{Q})$ in the Leray filtration.
 A multiplicative decomposition isomorphism as in (\ref{isodec}) induces by taking cohomology another decomposition
 of $H^2(\mathbb{P}(\mathcal{E}),\mathbb{Q})$ as $\pi^*H^2(B,\mathbb{Q})\oplus \mathbb{Q}h'$, where
 $h'=h+\pi^*\alpha$, for some $\alpha\in H^2(B,\mathbb{Q})$.
 In this multiplicative decomposition, $h'$ will generate a summand isomorphic to
 $H^0(B, R^2\pi_*\mathbb{Q})$. Let $r=rank\,\mathcal{E}$. As $c_1^{top}(\mathcal{E})=0$, one has $\pi_*h^r=0$ in $H^2(B,\mathbb{Q})$.
 As $(h')^r=0$ in $H^0(B, R^{2r}\pi_*\mathbb{Q})$, and $(h')^r$ belongs by multiplicativity to a direct summand
naturally isomorphic (by restriction to fibers) to $H^0(B, R^{2r}\pi_*\mathbb{Q})=0$,
 one must also have  $(h')^r=0$ in $H^{2r}(\mathbb{P}(\mathcal{E}),\mathbb{Q})$.
 On the other hand $(h')^r=h^r+r h^{r-1}\pi^*\alpha+\ldots +\pi^*\alpha^r$,
 and it follows that
 $$\pi_*(h')^r=0=\pi_*h^r+r \alpha\,\,{\rm in}\,\,H^2(B,\mathbb{Q}).$$
 Thus $\alpha=0$, $h'=h$, and $h^r=0$  in $H^{2r}(\mathbb{P}(\mathcal{E}),\mathbb{Q})$.
 The definition of Chern classes  and the fact that $h^r=0$
  show then that
$c_i^{top}(\mathcal{E})=0$ for all $i>0$.
\cqfd
In this example, the obstructions to the existence of a multiplicative
decomposition isomorphism are given by cycle classes on $B$. These classes vanish
locally on $B$ for the Zariski topology and this suggests studying the following variant
of Question \ref{question}:
\begin{question} \label{question2} Given a family of smooth projective varieties $\pi:\mathcal{X}\rightarrow B$, does there exist
a  Zariski dense open set $B^0$ of $B$, and a multiplicative decomposition isomorphism
as in (\ref{isodec}) for the restricted family $\mathcal{X}^0\rightarrow B^0$?
\end{question}
Our last example is given by families of curves and
shows that already in this case, we can have a negative answer to this weakened question. We fix an abelian surface,  choose a Lefschetz pencil of   curves
$C_t\subset A,\,t\in  \mathbb{P}^1$, and let  $B\subset \mathbb{P}^1$ be  the open set parameterizing
smooth fibers.
\begin{prop} \label{examplecurve}The family
$\pi:\mathcal{C}\rightarrow B$ does not admit  a multiplicative decomposition isomorphism over any non empty Zariski open set of  $B$.
\end{prop}
{\bf Proof.}
Assume there is a multiplicative decomposition isomorphism for the restricted family
$\pi:\mathcal{C}^0\rightarrow B^0$ over some non-empty Zariski open set $B^0$ of $B$. Then
we get by taking cohomology a decomposition
$$ H^1(\mathcal{C}^0,\mathbb{Q})\cong \pi^*H^1(B^0,\mathbb{Q})\oplus K,$$
where $K\cong H^0(B^0,R^1\pi_*\mathbb{Q})$ has the property
that
 the cup-product map:
$$\mu: K\otimes K \rightarrow H^2(\mathcal{C}^0,\mathbb{Q})$$
factors through the cup-product map
$$\mu:H^0(B,R^1\pi_*\mathbb{Q})\otimes H^0(B,R^1\pi_*\mathbb{Q})\rightarrow H^0(B,R^2\pi_*\mathbb{Q}).
$$

Now let $\alpha,\,\beta\in H^1(A,\mathbb{C})$ be the classes of two independent sections of $\Omega^1_A$.
Let us denote by $q:\mathcal{C}\rightarrow A$ the natural map.
Then we can decompose
$$q^*\alpha=\alpha_K+\pi^*\alpha',\,q^*\beta=\beta_K+\pi^*\beta',$$
with $\alpha_K,\,\beta_K\in K$ and $\alpha',\,\beta'\in H^1(B^0,\mathbb{C})$.
Taking their cup-product, and using  the fact that the cup-product
 is trivial on the summand $\pi^*H^1(B^0,\mathbb{C})$, we get the equality
$$q^*(\alpha\cup\beta)=\alpha_K\cup\beta_K+\alpha_K\cup\pi^*\beta'+\pi^*\alpha'\cup\beta_K,$$
and the first term $\alpha_K\cup\beta_K$ vanishes because it vanishes in
$H^0(B^0,R^2\pi_*\mathbb{C})$ (indeed, the classes $\alpha,\,\beta$ are of type $(1,0)$ and so are
their restrictions to the fibers $\mathcal{C}_b$ which are $1$-dimensional). The same arguments show that
$$q^*(\alpha\cup\beta)=q^*\alpha\cup \pi^*\beta' + \pi^*\alpha'\cup q^*\beta\,\,\,{\rm in}\,\,\,H^2(\mathcal{C}^0,\mathbb{C}).$$

The contradiction comes from the fact that $q^*(\alpha\cup\beta)$ does not vanish in
 $H^2(\mathcal{C}^0,\mathbb{C})$ (because this is the restriction of the class of a nonzero $(2,0)$-form on a projective completion of $\mathcal{C}^0$, namely the blow-up of $A$ at the base-points of the pencil)
and has trivial residues along all fibers $C_b,\,\,\,b\in \mathbb{P}^1\setminus B^0$, while
the independence of the restrictions of the classes $\alpha,\,\beta$ to the fibers $C_b,\,\,b\in \mathbb{P}^1\setminus B^0$
implies that the term on the right can have trivial residues along all fibers if and only if
$\beta'$ and $\alpha'$ have trivial residues at all points $b\in \mathbb{P}^1\setminus B^0$, which implies
$\beta'=0,\,\alpha'=0$.

\cqfd

Our main result in this paper is:
\begin{theo} \label{main14juilletintro}

(i)  For any smooth projective family  $\pi:\mathcal{X}\rightarrow B$
of  $K3$ surfaces,
there exist a  decomposition isomorphism
as in (\ref{isodec}) and a nonempty
 Zariski open subset $B^0$ of $B$, such that this decomposition becomes multiplicative
  for the restricted family $\pi:\mathcal{X}^0\rightarrow B^0$.

(ii) The class of the relative diagonal $[\Delta_{\mathcal{X}^0/B^0}]\in H^{4}(\mathcal{X}^0\times_{B^0}\mathcal{X}^0,\mathbb{Q})$ belongs
to the direct summand $H^0(B^0,R^4(\pi,\pi)_*\mathbb{Q})$ of $H^4(\mathcal{X}^0\times_{B^0}\mathcal{X}^0,\mathbb{Q})$, for the induced
decomposition of $R(\pi,\pi)_*\mathbb{Q}$.

(iii)  For any algebraic line bundle $\mathcal{L}$ on $\mathcal{X}$, there is a dense Zariski open set
$B^0$ of $B$ such that the topological Chern class $c_1^{top}(\mathcal{L})\in H^2(\mathcal{X},\mathbb{Q})$
restricted to $\mathcal{X} ^0$ belongs to the direct summand $H^0(B^0,R^2\pi_*\mathbb{Q})$
of $H^2(\mathcal{X}^0,\mathbb{Q})$ induced by this decomposition.
\end{theo}

 Statement  (i) is definitely wrong if we do not restrict to a Zariski open set (cf. section \ref{sec1.1} for an example).
Statement (iii) is in fact implied by (i), according to Lemma \ref{newlemma10aout}.

We note that statements (i) and (iii) together imply that the decomposition above coincide
locally over $B$ in the Zariski topology with the first one defined by Deligne in \cite{delignemotivedec}. This follows from
the characterization of the latter given in \cite[Prop. 2.7]{delignemotivedec}.

 We will explain in section \ref{sectionlink} how Theorem \ref{main14juilletintro} is related to the results of \cite{beauvoi}, \cite{beauville}
 (see also \cite{voisin} for further developments) concerning the Chow ring of $K3$ surfaces.
In fact the statement  was  motivated by the following
 result, which is an easy consequence of the results of
 \cite{beauvoi}, but can be seen as well as a consequence of Theorem \ref{main14juilletintro} by Proposition
 \ref{propproduitdecomp}.
 \begin{prop} \label{coromotiv} Let $\pi: \mathcal{S}\rightarrow B$ be a family of $K3$ surfaces,
   $\mathcal{L}_i\in {\rm Pic}\,\mathcal{S}$ and $n_{ij}$ be integers.
 Assume that the degree $4$ cohomology class $c=\sum_{ij}n_{ij}\,c_1^{top}(\mathcal{L}_i)c_1^{top}(\mathcal{L}_j)\in H^4(\mathcal{S},\mathbb{Q})$
 has trivial restriction on the fibers $\mathcal{S}_t,\,t\in B$ (or equivalently, has trivial restriction on one fiber $\mathcal{S}_t$, if
 $B$ is connected). Then there exists a nonempty
 Zariski open subset $B^0$ of $B$ such that $c$ vanishes in $H^4(\mathcal{S}^0,\mathbb{Q})$,
 where  $\mathcal{S}^0:=\pi^{-1}(B^0)$.
 \end{prop}

In section \ref{sectionlink}, we   prove    Proposition \ref{propproduitdecomp}, which says
in particular  that
Proposition \ref{coromotiv} is satisfied more generally by any family $\mathcal{X}\rightarrow B$ of varieties with trivial irregularity,  admitting a multiplicative decomposition
isomorphism, and for any fiberwise polynomial cohomological  relation between Chern classes of line bundles
on $\mathcal{X}$. This strongly relates the present work to the paper \cite{beauville}.

We will also use this proposition
  in section \ref{secnondec}  to
provide further examples of families of surfaces for which there is no multiplicative decomposition isomorphism
over any dense Zariski open set of the base, although there is no variation of Hodge structures in the fibers.

Let us mention one consequence of Theorem \ref{main14juilletintro}. Let $\pi:\mathcal{X}\rightarrow B$ be a projective family
of  $K3$ surfaces, with $B$ irreducible, and $\mathcal{L}\in {\rm Pic}\,\mathcal{X}$. Consider  the  $0$-cycle
$o_\mathcal{X}:=\frac{1}{{\rm deg}_{\mathcal{X}_t}\,\mathcal{L}^2}\mathcal{L}^2\in CH^2(\mathcal{X})_\mathbb{Q}$. By theorems \ref{BV} and
\ref{genprinciple}, this $0$-cycle
is independent of $\mathcal{L}$, at least after restriction to $\mathcal{X}^0=\pi^{-1}(B^0)$, for an adequate Zariski dense
open set $B^0$ of $B$. We also have the relative diagonal $\Delta_{\mathcal{X}/B}\in CH^2(\mathcal{X}\times_{B}\mathcal{X})$.
Let   $L_s$, $s\in I$, be line bundles on $\mathcal{X}$. Set $\mathcal{X}^{m/B}:=\mathcal{X}\times_B\ldots\times_B\mathcal{X}$, $\pi^{m}:\mathcal{X}^{m/B}\rightarrow B$, the $m$-th fibered product of $\mathcal{X}$ over $B$.
\begin{coro} \label{corepourhyper}  Consider a codimension $2r$ cycle $Z$ with $\mathbb{Q}$-coefficients in
$\mathcal{X}^{m/B}$ which is
a polynomial in the cycles $pr_i^*o_\mathcal{X}$, $pr_j^*L_s$, $pr_{kl}^*\Delta_{\mathcal{X}/B}$, where $1\leq i,j,k,l\leq m$. Assume that
the restriction of $Z$ to one (equivalently, any) fiber
$\mathcal{X}_t^m$ is cohomologous to $0$.
Then there exists a dense Zariski open set $B^0$ of $B$ such that
$Z$ is cohomologous to $0$ in $(\mathcal{X}^0)^{m/B}$.

\end{coro}

{\bf Proof.} Indeed, it follows from Theorem  \ref{main14juilletintro}, (iii)  that over a dense  Zariski open set $B^0$, the
classes $c_1^{top}(L_s)\in H^2(\mathcal{X}^0,\mathbb{Q})$ belong to the direct summand
$H^0(B^0,R^2\pi_*\mathbb{Q})$ of $H^2(\mathcal{X}^0,\mathbb{Q})$ induced by the multiplicative decomposition isomorphism
of Theorem \ref{main14juilletintro}. By multiplicativity, we have that the class $[o_\mathcal{X}]$
belongs to the direct summand
$H^0(B^0,R^4\pi_*\mathbb{Q})$ of $H^4(\mathcal{X}^0,\mathbb{Q})$. Theorem \ref{main14juilletintro}, (ii)
tells us that, over a Zariski open set $B^0$ of $B$, the class $[\Delta_{\mathcal{X}/B}]$ of the relative diagonal belongs
to the direct summand
$H^0(B^0,R^4(\pi,\pi)_*\mathbb{Q})$ of $H^4(\mathcal{X}^0\times_{B^0}\mathcal{X}^0,\mathbb{Q})$.
We thus conclude by multiplicativity that the class
$[Z]$ belongs to the direct summand
$H^0(B^0,R^{2r}(\pi^m)_*\mathbb{Q})$ of $H^{2r}(\mathcal{X}^0)^{m/B},\mathbb{Q})$.
But by assumption, the class $[Z]$ projects to $0$ in $H^0(B^0,R^{2r}(\pi^m)_*\mathbb{Q})$. We thus deduce that it is identically $0$.
\cqfd
 This corollary provides an evidence (of a rather speculative nature, in the same spirit as \cite{huy})
 for the conjecture made in \cite[Conj. 1.3]{voisin} concerning
 the Chow ring of hyper-K\"ahler manifolds, at least for those of type $S^{[n]}$, where $S$ is a $K3$ surface.
 Indeed, this conjecture states the following:
 \begin{conj}\label{conjrechau}  Let $Y$ be an algebraic hyper-K\"{a}hler variety. Then any polynomial cohomological relation
$P([c_1(L_j )], [c_i(T_Y )]) = 0\,\,{\rm  in}\,\, H^{2k}(Y,\mathbb{Q})$, $ L_j \in{\rm Pic}\, Y$,
already holds at the level of Chow groups :
$P(c_1(L_j), c_i(T_Y )) = 0 $ in $CH^k(Y )_\mathbb{Q}$.
 \end{conj}
Indeed, we proved in \cite[Prop. 2.5]{voisin} that for $Y=S^{[n]}$, this conjecture is implied by
the following conjecture:
 \begin{conj}\label{conjrechauchau} Let $S$ be an algebraic $K3$ surface. For any integer $m$, let
$P \in CH^p(S^m)_\mathbb{Q}$ be a weighted degree $k$ polynomial expression in
$pr_i^* c_1(L_s)$, $L_s\in {\rm  Pic}\, S$, $ pr_{jl}^*\Delta_S$:
Then if $[P] = 0$ in $H^{2k}(S^m,\mathbb{Q})$, we have $P = 0$ in $CH^k(S^m)_\mathbb{Q}$.
\end{conj}
 By
the general principle \ref{genprinciple}, Conjecture  \ref{conjrechauchau} implies Corollary \ref{corepourhyper}.
In the other direction, we can say the following (which is rather speculative): In the situation of
Conjecture  \ref{conjrechauchau}, we can find a family $\mathcal{X}\rightarrow B$
of smooth projective $K3$ surfaces, endowed with line
bundles $\mathcal{L}_s\in {\rm Pic}\,\mathcal{X}$, where everything is defined over $\mathbb{Q}$, such that
$S$ and the $L_s$'s are the fiber over some $t\in B$ of $\mathcal{X}$ and the $\mathcal{L}_j$'s.
Then we can construct using the same polynomial expression
the cycle
$\mathcal{P}\in CH^k(\mathcal{X}^{m/B})_\mathbb{Q}$
and Corollary \ref{corepourhyper} tells that the class of this cycle vanishes
in $H^{2k}((\mathcal{X}^0)^{m/B},\mathbb{Q})$. As $(\mathcal{X}^0)^{m/B}$ and $\mathcal{P}$
 are defined over $\mathbb{Q}$, the Beilinson conjecture predicts that it is trivial if furthermore its
 Abel-Jacobi invariant vanishes, which is presumably provable by the same method used to get the vanishing of the cycle class.

Theorem \ref{main14juilletintro} will be proved in section \ref{section2}. In fact, we will give two proofs of it.
In the first one, we use the existence of non trivial self $K$-correspondences
(see \cite{voisinKcorresp}), whose action on cohomology allows to split the cohomology in different pieces, in a way
which is compatible with the cup-product. This is very similar to the proof given in the abelian case (Proposition \ref{abelian}),
for which one uses homotheties.
The second proof is formal, and uses a curious decomposition of the small diagonal
$\Delta\subset S^3$
of a $K3$ surface $S$, obtained in \cite[Prop. 3.2]{beauvoi} (see  Theorem \ref{decompdiagonalpetiteintro}).

In section \ref{section3}, we will investigate the case of Calabi-Yau hypersurfaces $X$
in projective space $\mathbb{P}^n$ and establish for them
the following analogue of this decomposition of the small diagonal. We denote by $\Delta\cong X\subset X^3$ the small diagonal of $X$
and $\Delta_{ij}\cong X\times X\subset X^3$ the inverse image in $X^3$ of the diagonal of $X\times X$ by the
projection onto the product of the $i$-th and $j$-th factors. There is a natural $0$-cycle $o:=\frac{c_1(\mathcal{O}_X(1))^{n-1}}{n+1}\in CH_0(X)$.
\begin{theo} \label{theonewpourarticleintro} (cf. Theorem \ref{theonewpourarticle})
The following relation  holds in $CH^{2n-2}(X\times X\times X)_\mathbb{Q}$ (in the following equation,  ``$+(perm.)$'' means that we symmetrize in the indices the considered expression):
 \begin{eqnarray}\label{equadeltapetiteCYintro}\Delta=\Delta_{12}\cdot o_3+(perm.)+Z+\Gamma'\,\,{\rm in}\,\,CH^{2n-2}(X\times X\times X)_\mathbb{Q},
\end{eqnarray}
where $Z$ is the restriction to $X\times X\times X$
of a cycle of $\mathbb{P}^n\times\mathbb{P}^n\times\mathbb{P}^n$, and $\Gamma'$ is a multiple
of the following effective cycle of dimension $n-1$:
$$\Gamma:=\cup_{t\in F(X)}\mathbb{P}^1_t\times\mathbb{P}^1_t\times\mathbb{P}^1_t,$$
where $F(X)$ is the variety of lines contained in $X$.
 \end{theo}
As a consequence, we get the following result concerning the Chow ring of a Calabi-Yau hypersurface $X$ in $\mathbb{P}^n$, which generalizes
Theorem  1 of \cite{beauvoi} (see Theorem \ref{BV}):
 \begin{theo} Let $X$ be as above and let $Z_i,\,Z'_i$ be cycles of codimension $>0$ on $X$ such that ${\rm codim}\,Z_i+{\rm codim}\,Z'_i= n-1$. Then if we have a cohomological relation
 $$\sum_i n_i{\rm deg}\,(Z_i\cdot Z'_i)=0,$$
 this relation already holds at the level of Chow groups:
 $$\sum_i n_iZ_i\cdot Z'_i=0\,\,{\rm in}\,\,CH_0(X)_\mathbb{Q}.$$
 \end{theo}

 We conjecture that the cycle $\Gamma$ also comes from a cycle on $\mathbb{P}^n\times\mathbb{P}^n\times\mathbb{P}^n$. This would imply
 the analogue of Theorem \ref{main14juilletintro} for families of Calabi-Yau hypersurfaces.

\vspace{0.5cm}

{\bf Thanks.}  I thank Bernhard Keller for his help in the proof of Lemma \ref{ledecomp},  Christoph
Sorger and Bruno Kahn for useful discussions, and the  referee on a primitive version of this paper for useful comments.

\section{Link with the results of \cite{beauville}, \cite{beauvoi}\label{sectionlink}}
In this section, we first show how to deduce  Proposition  \ref{coromotiv}
      from the following Theorem
 proved in \cite{beauvoi} :

 \begin{theo} \label{BV} (Beauville-Voisin 2004)  Let $S$ be a $K3$ surface, $D_i\in CH^1(S)$ be divisors on $S$
 and $n_{ij}$ be integers. Then if the $0$-cycle $\sum_{i,j}n_{ij}D_i D_j\in CH_0(S)$ is cohomologous to $0$ on $S$, it is equal to $0$ in $CH_0(S)$.

 \end{theo}

We will use here and many times later on in the paper the following ``general principle'' (cf.
\cite{blochsrinivas}, \cite[Theorem 10.19]{voisinbook}, \cite[Corollary 3.1.6]{voisinweyllectures}:
\begin{theo}\label{genprinciple} Let $\pi:X\rightarrow B$  be
 a morphism with $X,\,B$ smooth, and $Z\in CH^k(X)$ such that
 $Z_{\mid X_t}=0$ in $CH^k(X_t)$
 for any $t\in B$. Then there exists a dense Zariski open set $B^0\subset B$ such that
 \begin{eqnarray}
 \label{cirquecoro}
  [Z]  =0\,\,{\rm in}\,\,{H}^{2k}(X^0,\mathbb{Q}),\end{eqnarray}
  where $X^0:=\pi^{-1}(B^0)$.
\end{theo}

 {\bf  Proof of Proposition  \ref{coromotiv}.}
 Indeed,
  under the  assumption that the intersection number $\sum_{i,j} n_{ij}c_1^{top}(\mathcal{L}_{i,b})c_1^{top}(\mathcal{L}_{j,b})=0$
 vanishes in  $H^4(\mathcal{S}_b,\mathbb{Q})=\mathbb{Q}$  for all $b\in B$, Theorem \ref{BV} says that the codimension
 $2$ cycle $\sum_{i,j}n_{ij}c_1(\mathcal{L}_{i})c_1(\mathcal{L}_j)\in CH^2(\mathcal{S})$ has trivial restriction
 on each fiber $\mathcal{S}_b$. The general principle
 \ref{genprinciple}
 then implies that there is a Zariski dense open set $B^0$ of $B$ such that
  the class $\sum_{i,j}n_{ij}c_1^{top}(\mathcal{L}_{i})c_1^{top}(\mathcal{L}_j)$ vanishes in
   $H^4(\mathcal{S}^0,\mathbb{Q})$.
\cqfd

We next prove the following Proposition \ref{propproduitdecomp}, which
provides a conclusion similar as above, under the assumption that the family has
a multiplicative decomposition isomorphism over a Zariski open set.

   Let  $\pi: \mathcal{X}\rightarrow B$ be a projective family of smooth  complex  varieties
   such that  $H^1(\mathcal{X}_b,\mathcal{O}_{\mathcal{X}_b})=0$ for any $b\in B$,
 parameterized by a connected complex quasi-projective  variety $B$.  Let $\mathcal{L}_i,\,i=1,\ldots,m$ be line bundles on
  $\mathcal{X}$ and $l_i:=c_1^{top}(\mathcal{L}_i)\in H^2(\mathcal{X},\mathbb{Q})$.
  We will say that a cohomology class $\beta\in H^*(\mathcal{X},\mathbb{Q})$ is Zariski locally trivial over
  $B$ if $B$ is covered
 by Zariski open sets
  $B^0\subset B$, such that
 $\beta_{\mid \mathcal{X}^0}=0$ in $H^{*}(\mathcal{X}^0,\mathbb{Q})$, where
  $\mathcal{X}^0=\pi^{-1}(B^0)$.
  \begin{prop} \label{propproduitdecomp} Assume that there is a multiplicative decomposition isomorphism
 \begin{eqnarray}\label{9aout}R\pi_*\mathbb{Q}=\oplus_iR^i\pi_*\mathbb{Q}[-i].
 \end{eqnarray}
 Let $P$ be a homogeneous polynomial of degree $r$ in $m$ variables with rational coefficients and let
 $\alpha:=P(l_i)\in H^{2r}(\mathcal{X},\mathbb{Q})$.
 Then, if
 $\alpha_{\mid \mathcal{X}_b}=0\,\,{\rm in}\,\, H^{2r}(\mathcal{X}_b,\mathbb{Q})$ for some  $b\in B$,
 the class  $\alpha$ is  Zariski locally trivial over
  $B$.

 \end{prop}

 {\bf Proof.}
  We will assume for simplicity that $B$ is smooth although a closer look at the proof shows that this assumption is not necessary.
 The multiplicative decomposition isomorphism induces, by taking cohomology and using the fact that
 the fibers have no degree $1$ rational cohomology, a decomposition
\begin{eqnarray}\label{27noveqn0} H^2(\mathcal{X},\mathbb{Q})= H^0(B,R^2\pi_*\mathbb{Q})\oplus \pi^*H^2(B,\mathbb{Q}),
\end{eqnarray}
 which is compatible with cup-product, so that the cup-product
 map on the first term factors through the map induced by cup-product:
 $$\mu_r: H^0(B,R^2\pi_*\mathbb{Q})^{\otimes r}\rightarrow H^0(B,R^{2r}\pi_*\mathbb{Q}).$$
 We write  in this decomposition $l_i=l'_i+\pi^*k_i$, where $k_i\in H^2(B,R^0\pi_*\mathbb{Q})=H^2(B,\mathbb{Q})\stackrel{\pi^*}{\cong}\pi^*H^2(B,\mathbb{Q})$.
 We now have:
 \begin{lemm} \label{newlemma10aout}
 The assumptions being as in  Proposition \ref{propproduitdecomp}, the classes $k_i$ are divisor classes on $B$. Thus $B$
 is covered by Zariski open sets $B^0$ such that the
 divisor classes $l_i$ restricted to $\mathcal{X}^0$ belong  to the direct summand $H^0(B^0,R^2\pi_*\mathbb{Q})$.
 \end{lemm}
  {\bf Proof.}
  Indeed, take any line bundle  $\mathcal{L}$
  on $\mathcal{X}$.  Let $l=c_1^{top}(\mathcal{L})\in H^2(\mathcal{X},\mathbb{Q})$ and decompose
  as above $l=l'+\pi^*k$, where $l'$ has the same image as $l$ in $H^0(B,R^2\pi_*\mathbb{Q})$ and $k$
  belongs to $ H^2(B,\mathbb{Q})$.
  Denoting by  $n$  the dimension of the fibers, we get:
  \begin{eqnarray}\label{27noveqn}
   l^nl_i=(\sum_{p}{n\choose p}{l'}^p\pi^*k^{n-p})(l'_i+\pi^*k_i)=\sum_{p}{n\choose p}{l'}^pl'_i\pi^*k^{n-p}+\sum_{p}{n\choose p}{l'}^p\pi^*(k^{n-p}k_i).
   \end{eqnarray}
  Recall now that the decomposition is multiplicative. The class ${l'}^nl_i'$ thus belongs to the direct summand
  of $H^{2n+2}(\mathcal{X},\mathbb{Q})$ isomorphic to $H^0(B,R^{2n+2}\pi_*\mathbb{Q})$  deduced from the decomposition
  (\ref{9aout}). As $R^{2n+2}\pi_*\mathbb{Q}=0$, we conclude that ${l'}^nl_i'=0$.
   Applying $\pi_*: H^{2n+2}(\mathcal{X},\mathbb{Q})\rightarrow H^2(B,\mathbb{Q})$ to (\ref{27noveqn}), we then get:
  \begin{eqnarray}\label{27noveqn2}
  \pi_*(l^nl_i)=n{\rm deg}_{\mathcal{X}_b}({l'}^{n-1}l'_i)k+{\rm deg}_{\mathcal{X}_b}({l'}^n)k_i=n{\rm deg}_{\mathcal{X}_b}({l}^{n-1}l_i)k+{\rm deg}_{\mathcal{X}_b}({l}^n)k_i.
  \end{eqnarray}
  Observe that the term on the left is a divisor class on $B$.
  If the fiberwise self-intersection ${\rm deg}_{\mathcal{X}_b}({l_i}^n)$ is non zero, we can take $\mathcal{L}=\mathcal{L}_i$
  and (\ref{27noveqn2}) gives:
  $$\pi_*(l_i^{n+1})=(n+1){\rm deg}_{\mathcal{X}_b}({l_i}^n)k_i.$$
  This shows that $k_i$ is a divisor class on $B$ and proves the lemma in this case.
  If ${\rm deg}_{\mathcal{X}_b}({l_i}^n)$ is equal to  $0$, choose a line bundle
  $\mathcal{L}$ on $\mathcal{X}$ such that both intersection numbers
   ${\rm deg}_{\mathcal{X}_b}({l}^{n-1}l_i)$ and ${\rm deg}_{\mathcal{X}_b}({l}^n)$ are nonzero
  (such an $\mathcal{L}$ exists because the morphism
  $\pi$ is projective). Then, in the formula
  $$\pi_*(l^nl_i)= n{\rm deg}_{\mathcal{X}_b}({l}^{n-1}l_i)k+{\rm deg}_{\mathcal{X}_b}({l}^n)k_i,$$
  the left hand side is a divisor class on $B$ and, as we just proved, the first term in the right hand side is also a divisor
  class on $B$. It thus follows that ${\rm deg}_{\mathcal{X}_b}({l}^n)k_i$ is
   a divisor class on $B$. The lemma
   is thus proved.
\cqfd

Coming back to the proof of Proposition \ref{propproduitdecomp},
Lemma \ref{newlemma10aout} tells us that   $B$ is covered by Zariski open sets $B^0$ on which
$l_i$ belongs   to the first summand $H^0(B^0,R^2\pi_*\mathbb{Q})$ in (\ref{27noveqn0}).  It then follows
by multiplicativity that any polynomial expression $P(l_i)_{\mid \mathcal{X}^0}$ belongs to a direct summand
of $H^{2r}(\mathcal{X}^0,\mathbb{Q})$ isomorphic by the natural projection to $H^0(B^0,R^{2r}\pi_*\mathbb{Q})$.
Consider now our fiberwise cohomological polynomial relation $\alpha_{\mid \mathcal{X}_b}=0$ in $H^{2r}(\mathcal{X}_b,\mathbb{Q})$, for some
 $b\in B$. Since $B$ is connected, it says equivalently that  $\alpha$ vanishes in
 $H^0(B^0,R^{2r}\pi_*\mathbb{Q})$. It follows then from the previous statement that it vanishes in  $H^{2r}(\mathcal{X}^0,\mathbb{Q})$.
 \cqfd
\subsection{Application\label{secnondec}}
We can use Proposition \ref{propproduitdecomp} to exhibit very simple
families of smooth projective surfaces, with no variation of Hodge structure, but for which
there is no multiplicative decomposition isomorphism on any nonempty Zariski open set of the base.

We consider  a smooth projective surface $S$, and
set
$$X=\widetilde{(S\times S)}_\Delta,\,B=S,\,\pi=pr_2\circ \tau,$$
where
$\tau:\widetilde{(S\times S)}_\Delta\rightarrow S\times S$ is the blow-up of the diagonal.
\begin{prop} \label{examplediag}
Assume that $h^{1,0}(S)=0,\,h^{2,0}(S)\not=0$. Then there is no
multiplicative decomposition isomorphism for $\pi:X\rightarrow B$ over any Zariski dense open set of $B=S$.
\end{prop}
{\bf Proof.} Let $H$ be an ample line bundle on $S$, and $d:={\rm deg}\,c_1(H)^2$.
On $X$, we have then two line bundles, namely
$L:=\tau^*(pr_1^*H)$ and $L'=\mathcal{O}_X(E)$ where $E$ is the exceptional divisor of $\tau$.
On the fibers of $\pi$, we have the relation
$$ {\rm deg}\,c_1(L)^2=-d\, {\rm deg}\,c_1(L')^2.$$
If there existed a
multiplicative decomposition isomorphism over a Zariski dense open set of $B=S$,
we would have by Proposition \ref{propproduitdecomp}, using the fact that
the fibers of $\pi$ are regular, a Zariski dense open set $U\subset S$ such that
the relation
\begin{eqnarray}\label{rel12juill} c_1^{top}(L)^2=-d \,c_1^{top}(L')^2
\end{eqnarray}
holds in $H^4(X_U,\mathbb{Q})$.
If we apply $\tau_*:H^4(X_U,\mathbb{Q})\rightarrow H^4(S\times U,\mathbb{Q})$ to this relation, we now get:
\begin{eqnarray}\label{rel12juilldeux} pr_1^*c_1^{top}(H)^2=-d [\Delta]
\end{eqnarray}
in $H^4(S\times U,\mathbb{Q})$.

This relation implies that the class
$pr_1^*c_1^{top}(H)^2+d [\Delta]\in H^4(S\times S,\mathbb{Q})$ comes
from a class  $\gamma\in H^2(
S\times \widetilde{D},\mathbb{Q})$, where $D:=S\setminus U$
and $\widetilde{D}$ is a desingularization of $D$.
 Denoting by $\tilde{j}:\widetilde{D}\rightarrow S$
 the natural map, we then conclude that for any class $\alpha\in H^2(S,\mathbb{Q})$,
 $$d\alpha\in H^2(S,\mathbb{Q})=\tilde{j}_*(\gamma_*\alpha)$$ is supported
 on $D$. This contradicts  the assumption $h^{2,0}(S)\not=0$.
\cqfd

 \subsection{Example where Theorem \ref{main14juilletintro}, (i) is not satisfied globally on $B$ \label{sec1.1}}
 Let us apply the same arguments as in the proof of Proposition
 \ref{propproduitdecomp} to exhibit simple
 families of smooth projective $K3$ surfaces for which a multiplicative decomposition isomorphism
 does not exist on the whole base.

 We take $B=\mathbb{P}^1$ and $\mathcal{S}\subset \mathbb{P}^1\times \mathbb{P}^1\times \mathbb{P}^2$
 a generic hypersurface of multidegree $(d,2,3)$.
 We put $\pi:=pr_1$. This is not a smooth family of $K3$ surfaces because of the nodal fibers, but we can take a
   finite cover of $\mathbb{P}^1$ and introduce a simultaneous resolution  of the pulled-back family to get a  family of smooth $K3 $ surfaces parameterized by a complete curve. (Note that
 the simultaneous resolution does not hold in the projective category, so the morphism
 $\pi':\mathcal{S} '\rightarrow B'$ obtained this way is usually not projective : this is a minor point.)
 By the Grothendieck-Lefschetz theorem, one has
 ${\rm NS}\,\mathcal{S}=\mathbb{Z}^3={\rm NS}\,(\mathbb{P}^1\times \mathbb{P}^1\times \mathbb{P}^2)$, and
 the same is true for $\mathcal{S}'$ if one assumes that the general fiber of $\mathcal{S}$ over $B$ has
 ${\rm Pic}\,\mathcal{S}_b={\rm Pic}\,(\mathbb{P}^1\times \mathbb{P}^2)=\mathbb{Z}^2$.

 We prove now :
 \begin{lemm} The family $\pi':\mathcal{S}'\rightarrow B'$ does not admit a multiplicative decomposition isomorphism over $B'$.
 \end{lemm}
 {\bf Proof.} As the hypersurface $\mathcal{S}$ is generic, the family $\mathcal{S}\rightarrow \mathbb{P}^1$
 is not locally isotrivial. It follows that $H^2(\mathcal{S}',\mathcal{O}_{\mathcal{S}'})=0$, and thus
 $$H^2(\mathcal{S}',\mathbb{Q})=NS(\mathcal{S}')\otimes\mathbb{Q}^.$$
 As already mentioned, the right hand side is isomorphic to $\mathbb{Q}^3$, generated
 by the pull-back to $\mathcal{S}'$ of the natural classes $h_1,\,h_2,\,h_3$ on
 ${\rm Pic}\,(\mathbb{P}^1\times \mathbb{P}^1\times \mathbb{P}^2)$. The first class
 $h_1$ belongs to the natural summand ${\pi'}^{*}H^2(B,\mathbb{Q})=H^2(B,R^0\pi'_*\mathbb{Q})$ and, as explained above, the existence of a multiplicative decomposition
 isomorphism would imply the existence of a decomposition
 $${\rm NS}\,(\mathcal{S}')\otimes\mathbb{Q}=H^2(\mathcal{S}',\mathbb{Q})
 ={\pi'}^*H^2(B',\mathbb{Q})\oplus H,\,H\cong H^0(B,R^2\pi'_*\mathbb{Q})$$
 such that the cup-product map on
 $H$ factors through the map given by cup-product
 $$\mu: H^0(B',R^2\pi'_*\mathbb{Q})\otimes  H^0(B',R^2\pi'_*\mathbb{Q})\rightarrow  H^0(B',R^4\pi'_*\mathbb{Q})=\mathbb{Q}.$$
 Let us show that such a decomposition does not exist. As all classes are obtained by pull-back from
 $\mathcal{S}$ and the pull-back map preserves the cup-product, we can make the computation on $\mathcal{S}$.
 Let $h'_2=h_2-\alpha h_1,\,h'_3=h_3-\beta h_1$ be generators for $H$.
 The class $h'_2$ has  self-intersection $0$ on the fibers $\mathcal{S}_b$, and
 it follows that we must have ${h'}_2^2=0$ in $H^4(\mathcal{S},\mathbb{Q})$.
 As $h_2^2=0$ and ${h'}_2^2=h_2^2-2\alpha h_1h_2$, with $h_1h_2\not=0$ in $H^4(\mathcal{S},\mathbb{Q})$,
 we conclude that $\alpha=0$ and $h_2=h'_2$.
 Next, the class $h_3^2$ (hence also
 the class ${h'}_3^2$) has degree $2$ on the fibers $\mathcal{S}_b$ ; furthermore the intersection number
 $h_2h_3$ of the classes $h_2$ and $h_3$ on the fibers $\mathcal{S}_b$ is equal to $3$  (thus we get as well
 that the intersection number
 $h_2h'_3=h'_2h'_3$ on the fibers $\mathcal{S}_b$ is equal to $3$).

 If our multiplicative decomposition exists, we conclude that we must have the following relation
 in $H^4(\mathcal{S},\mathbb{Q})$:
 \begin{eqnarray}\label{rel} 3 {h'}_3^2-2h_2h'_3=0.
 \end{eqnarray}
 Equivalently, as the class of $\mathcal{S}$ in $\mathbb{P}^1\times \mathbb{P}^1\times \mathbb{P}^2$
 is an ample class equal to $dh_1+2h_2+3h_3$, we should have:
  \begin{eqnarray}\label{rel2}(dh_1+2h_2+3h_3)(3 (h_3^2-2\beta h_3h_1)-2h_2(h_3-\beta h_1))=0\,\,{\rm in}\,\,H^6(\mathbb{P}^1\times \mathbb{P}^1\times \mathbb{P}^2,\mathbb{Q}).
  \end{eqnarray}
  However this class is equal to
  $(3d-18\beta)h_1h_3^2+(-2d-6\beta)h_1h_2h_3$, where the two classes
  $h_1h_3^2,\, h_1h_2h_3$ are independent in
  $H^6(\mathbb{P}^1\times \mathbb{P}^1\times \mathbb{P}^2,\mathbb{Q})$.
  We conclude that for the  equation (\ref{rel2}) to hold, one needs
  $$3d-18\beta=0,\,-2d-6\beta=0,$$
  which has no  solution for $d\not=0$.
  Hence the relation (\ref{rel}) is not satisfied for any choice of $h'_3$.

 \cqfd

 \section{Proof of Theorem \ref{main14juilletintro}\label{section2}}
 \subsection{A criterion for the existence of a decomposition}
Our proofs will be based on the following easy and presumably standard
lemma, applied to the category of sheaves of $\mathbb{Q}$-vector spaces on $B$.

Let ${A}$ be a $\mathbb{Q}$-linear abelian category, and let $\mathcal{D}({A})$ be the corresponding derived category
of left  bounded  complexes.
Let $M\in \mathcal{D}(A)$ be an object with bounded cohomology such that $End\,M$ is finite dimensional.
Assume $M$ admits a morphism $\phi: M\rightarrow M$ such that
$$H^i(\phi):H^i(M)\rightarrow H^i(M)$$
is equal to $\lambda_i Id_{H^i(M)}$, where all the $\lambda_i\in\mathbb{Q}$ are distinct.
\begin{lemm}  \label{ledecomp} The morphism $\phi$ induces a canonical decomposition
\begin{eqnarray}
\label{splittingle}
 M\cong \oplus_i H^i(M)[-i],
 \end{eqnarray}
  characterized by the properties :

  1) The induced map on cohomology is the identity map.

  2)  One has
   \begin{eqnarray}
 \label{chareigenvalue} \phi\circ \pi_i=\lambda_i \pi_i:\,M\rightarrow M.
  \end{eqnarray}
   where $\pi_i$ corresponds via the isomorphism (\ref{splittingle}) to the $i$-th projector $pr_i$.
\end{lemm}
{\bf Proof.} We first prove using the arguments of \cite{deligne} that $M$ is decomposed, namely there is an isomorphism
$$f:M\cong \oplus_i H^i(M)[-i].$$
For this, given an object $K\in Ob\,A$, we consider the left exact functor $T$ from ${A}$ to the category of $\mathbb{Q}$-vector spaces defined by
$T(N)= Hom_A(K, N)$, and for any integer $i$ the induced functor,  denoted by $T_i$,
 $N\mapsto Hom_{\mathcal{D}(A)}( K[-i], N)$ on $\mathcal{D}(A)$.
For any $N\in \mathcal{D}(A)$, there is the hypercohomology spectral sequence with $E_2$-term
$$E_2^{p,q}=R^pT_i(H^q(N))=Ext^{p+i}_A(K,H^q(N))\Rightarrow \mathbb{R}^{p+q}T_i(N).$$
Under our assumptions, this spectral sequence for $N=M$ degenerates at $E_2$.
Indeed, the morphism $\phi$ acts then on the above spectral sequence starting from $E_2$.
The differential
$d_2:E_2^{p,q}\rightarrow E_2^{p+2,q-1} $
\begin{eqnarray}
\label{d2expl} Ext^{p+i}_A(K,H^q(M))\Rightarrow Ext^{p+2+i}_A(K,H^{q-1}(M))
\end{eqnarray}
 commutes with the action of $\phi$.
 On the other hand, $\phi$ acts as $\lambda_q Id$ on the left hand side and as
$\lambda_{q-1} Id$ on the right hand side of (\ref{d2expl}). Thus we conclude that $d_2=0$ and similarly that all
$d_r,\,r\geq 2$ are $0$.

We take now $K=H^i(M)$. We conclude from the degeneracy at $E_2$ of the above spectral sequence that the map
$$ Hom_{\mathcal{D}(A)}\,(H^i(M)[-i],M)\rightarrow Hom_A(H^i(M), H^i(M))=E_2^{-i,i}$$
is surjective, so that there is a morphism
$$f_i: H^i(M)[-i]\rightarrow M$$
inducing the identity on degree $i$ cohomology. The  direct sum   $f=\sum f_i$  is a quasi-isomorphism which gives the desired splitting.

The morphism $\phi$ can thus be seen as a morphism of the split object
$\oplus_iH^i(M)[-i]$. Such a morphism is given by a block-uppertriangular  matrix
$$\phi_{j,i}\in Ext^{i-j}_A(H^i (M), H^j(M)),\,i\geq j,$$
with $\lambda_i Id$ on the $i$-th diagonal block.
Let $\psi$ be the endomorphism of $End\,M$ given by left multiplication by
$\phi$.
 We have by the above description of $\phi$:
  \begin{eqnarray} \label{characteristic}\prod_{i,H^i(M)\not=0}(\psi-\lambda_i Id_{End\,M})=0,
  \end{eqnarray}
   which shows that the endomorphism $\psi$ is diagonalizable. More precisely, as $\psi$ is block-uppertriangular in
    an adequately  ordered  decomposition
    $$End\,M=\oplus_{i\geq j}Ext^{i-j}_A(H^i (M), H^j(M)),$$ with  term $\lambda_j Id$ on the block diagonals
    $Ext^{i-j}_A(H^i (M), H^j(M))$, hence in particular on  $End_A\,H^j(M)$, we conclude that
    there exists $\pi'_i\in End\,M$ such that $\pi'_i$ acts as the identity on
    $H^i(M)$, and $\phi\circ \pi'_i=\lambda_i \pi'_i$.

Let $\rho_i:=\pi'_i\circ f_i: H^i(M)[-i]\rightarrow M$. Then
$\rho:=\sum\,\rho_i$ gives another decomposition
$\oplus_iH^i(M)[-i]\cong M$
and we have $\phi\circ\rho_i=\lambda_i\rho_i$, which gives
$\phi\circ \pi_i=\lambda_i\pi_i$, where $\pi_i=\rho\circ pr_i\circ \rho^{-1}$.

 The uniqueness of the $\pi_i$'s satisfying properties 1) and 2) is obvious, since
 these properties force the equality
 $\pi_i=\frac{\prod_{j\not=i}(\phi-\lambda_j Id_M)}{\prod_{j\not=i}\lambda_i-\lambda_j}.$

\cqfd
The following result is proved in \cite{deningermurre}, by similar but somehow more complicated arguments (indeed they use Fourier-Mukai transforms, which exist only in the projective case):
\begin{coro} (Deninger-Murre 1991) \label{corab} For any family $\pi:\mathcal{A}\rightarrow B$ of abelian varieties or complex tori, there is
a multiplicative decomposition isomorphism $R\pi_*\mathbb{Q}=\oplus_i R^i\pi_*\mathbb{Q}[-i]$.
\end{coro}
{\bf Proof.} Choose an integer $n\not=\pm 1$ and consider
the multiplication map
$$\mu_n:\mathcal{A}\rightarrow \mathcal{A},\,\,\, a\mapsto na.$$
We then get morphisms
$\mu_n^*:R\pi_*\mathbb{Q}\rightarrow R\pi_*\mathbb{Q}$
with the property that the induced morphisms on each $R^i\pi_*\mathbb{Q}=H^i(R\pi_*\mathbb{Q})$ is multiplication
by $n^i$.
We use now Lemma \ref{ledecomp}  to deduce from
 such a morphism  a canonical splitting \begin{eqnarray}
 \label{splitting}
 R\pi_*\mathbb{Q}\cong \oplus_iR^i\pi_*\mathbb{Q}[-i],
 \end{eqnarray}
  characterized by the properties that the induced map on cohomology is the identity map, and
   \begin{eqnarray}
 \label{eigenvalue} \mu_n^*\circ \pi_i=n^i \pi_i:\,R\pi_*\mathbb{Q}\rightarrow R\pi_*\mathbb{Q}.
  \end{eqnarray}
   where $\pi_i$ is the endomorphism of $R\pi_*\mathbb{Q}$ which identifies to
   the $i$-th projector  via the isomorphism (\ref{splitting}).
On the other hand, the morphism  $ \mu : R\pi_*\mathbb{Q}\otimes R\pi_*\mathbb{Q}\rightarrow R\pi_*\mathbb{Q}$
given by cup-product
is compatible with $\mu_n^*$, in the sense that
$$\mu\circ (\mu_n^*\otimes\mu_n^*)=\mu_n^*\circ \mu:R\pi_*\mathbb{Q}\otimes R\pi_*\mathbb{Q}\rightarrow R\pi_*\mathbb{Q}.$$
Combining this last equation with  (\ref{eigenvalue}), we find that
$$\mu\circ (\mu_n^*\otimes\mu_n^*)\circ (\pi_i\otimes\pi_j)=n^{i+j}\mu\circ(\pi_i\otimes\pi_j)$$
$$=\mu_n^*\circ \mu\circ(\pi_i\otimes\pi_j):R\pi_*\mathbb{Q}\otimes R\pi_*\mathbb{Q}\rightarrow R\pi_*\mathbb{Q}$$
from which it follows applying again (\ref{eigenvalue}) that
$\mu\circ \pi_i\otimes\pi_j$ factors through  $R^{i+j}\pi_*[-i-j]$, or equivalently that
in  the splitting (\ref{splitting}), the cup-product morphism $\mu$
maps $ R^i\pi_*\mathbb{Q}[-i]\otimes R^j\pi_*\mathbb{Q}[-j])$ to the
summand $ R^{i+j}\pi_*[-i-j]$.
\cqfd
\subsection{$K$-autocorrespondences}

$K$-correspondences were introduced in \cite{voisinKcorresp} in order to study intrinsic volume forms on complex manifolds.
\begin{Defi} \label{defi}(Voisin 2004) A $K$-isocorrespondence between two projective complex manifolds $X$ and $Y$ of dimension $n$
 is a
 $n$-dimensional closed algebraic subvariety $\Sigma\subset X\times Y$,
 such
 that  each irreducible component of $\Sigma$ dominates $X$ and $Y$ by the natural  projections, and satisfying the
 following  condition :
Let $\widetilde\Sigma\stackrel{\tau}{\rightarrow}\Sigma$ be
 a desingularization, and let
 $f:=pr_1\circ\tau:\widetilde\Sigma\rightarrow X,\,g:=pr_2\circ
 \tau:\widetilde\Sigma\rightarrow Y$.
 Then we have the equality
 \begin{eqnarray}
 \label{eqram}R_f= R_g
 \end{eqnarray} of the  ramification divisors of $f$ and $g$ on
 $\widetilde\Sigma$.

A $K$-autocorrespondence of $X$ is a  $K$-isocorrespondence between  $X$ and itself.

\end{Defi}
We will be interested in $K$-autocorrespondences $\Sigma\subset X\times X$, where $X$ is a smooth
complex projective variety with trivial canonical bundle. In fact, we are not  interested
in this paper in the equality (\ref{eqram}) of ramification divisors, but in the
proportionality of pulled-back top holomorphic forms, which is an equivalent property  by the following lemma:
\begin{lemm} \label{leeqdef} Let $X$ be a smooth complex compact manifold with trivial canonical bundle, and
let $\Sigma\subset X\times X$
be an irreducible self-correspondence, with desingularization $\tau:\widetilde\Sigma\rightarrow \Sigma$. Then $\Sigma$ is a
$K$-autocorrespondence if and only if for some coefficient $\lambda\in\mathbb{C}^*$, one has

\begin{eqnarray}\label{propform} 0\not= f^*\eta=\lambda g^*\eta\,\,{\rm in}\,\,H^0(\widetilde\Sigma,K_{\widetilde\Sigma})
\end{eqnarray}
for any nonzero holomorphic section $\eta$ of $K_X$,
where as before $f=pr_1\circ\tau,\,g=pr_2\circ \tau$.
\end{lemm}
{\bf Proof.} Indeed, as $f^*\eta$ and $g^*\eta$ are not identically $0$,
the maps $f$ and $g$ are dominating and thus generically finite. As $K_X$ is trivial, $R_f$ and $R_g$ are respectively the divisors
of the  pulled-back forms $f^*\eta,\,g^*\eta\in H^0(\widetilde\Sigma,K_{\widetilde\Sigma})$.
As $\widetilde{\Sigma}$ is irreducible, these two forms are thus proportional if and only if $R_f=R_g$.
\cqfd
The simplest way to construct such a $K$-autocorrespondence is by studying rational equivalence of points on $X$: We recall for the convenience of the reader the proof of the following statement,
which can be found in \cite[Sec. 2]{voisinKcorresp}:
 Let $X$ be a complex projective $n$-fold with trivial canonical
bundle, and $z_0\in CH_0(X)$ be a fixed $0$-cycle. Let $m_1,\,m_2$ be non zero integers.
\begin{prop} \label{criKcor}Let $\Sigma\subset X\times X$ be a $n$-dimensional subvariety which dominates $X$ by both projections, and such that, for any $(x,y)\in \Sigma$, $m_1x+m_2y=z_0$ in $CH_0(X)$. Then $\Sigma$ is a $K$-autocorrespondence of
$X$. More precisely, we have the equality
$m_1f^*\eta=-m_2g^*\eta$ in $H^0(\widetilde{\Sigma},K_{\widetilde{\Sigma}})$ for any  holomorphic
$n$-form $\eta$ on $X$.
\end{prop}
{\bf Proof.} Let
$\tau:\widetilde{\Sigma}\rightarrow \Sigma$ be a desingularization of
$\Sigma$ and let as above $f:=pr_1\circ\tau,\,g=pr_2\circ
 \tau$.
  We apply Mumford's theorem \cite{mumford} or its generalization \cite[Proposition 10.24]{voisinbook} to the
  cycle
  $$\Gamma=m_1{\rm Graph}(f)+m_2{\rm Graph}(g)\in CH^n(\widetilde{\Sigma}\times X)$$ which has the property
   that ${\rm Im}(\Gamma_*:CH_0(\widetilde{\Sigma})_{hom}\rightarrow CH_0(X))$ is supported on  ${\rm Supp}\,z_0$.
  It follows that
for any holomorphic form $\eta$ of degree $>0$ on $X$,
$\Gamma^*\eta=0$ on $\widetilde{\Sigma}$. But  we have
$$\Gamma^*\eta=m_1f^*\eta+m_2g^*\eta \,\,{\rm in}\,\, H^0(\widetilde{\Sigma},\Omega_{\widetilde{\Sigma}}^l).$$
For  $l=n$, we get the desired
 equality $m_1f^*\eta=-m_2g^*\eta$ in $H^0(\widetilde{\Sigma},K_{\widetilde{\Sigma}})$.

\cqfd

Let $S$ be an algebraic  $K3$ surface, and $L$ an ample  line bundle on $S$ of self-intersection  $c_1(L)^2=2d$.
 We assume that ${\rm Pic}\,S$ has rank $1$, generated by
 a class proportional to $L$.
 There is a $1$-dimensional family of singular elliptic
curves in $|L|$  which sweep-out $S$. They may be not irreducible, and have in particular
fixed rational components, but as $({\rm Pic}\,S)\otimes\mathbb{Q}$ is generated by $L$, the classes of
all irreducible  components are proportional to $c_1(L)$. Changing $L$ if necessary, we may then assume the
general fibers of
  this $1$-dimensional family of elliptic curves are irreducible.
Starting from this one dimensional  family of irreducible elliptic curves $\Sigma_1:=\bigcup_{b\in \Gamma_1}{\Sigma'_b}$, we get by desingularizing
  $\Sigma_1$ and $\Gamma_1$ the following data:
 A smooth projective surface $\Sigma$, and  two morphisms
$$\phi:\Sigma \rightarrow S,\,p:\Sigma\rightarrow \Gamma,$$
 where $p$ is  surjective  with
elliptic fibers $\Sigma_b$  such that  $\phi_*(\Sigma_b)\in|L|$, $\Gamma$ is
a smooth curve, and $\phi$ is generically finite.

Choose an integer $m\equiv 1$  mod. $2d$, and write $m=2kd+1$. For a general point  $x\in \Sigma$, the fiber $\Sigma_x:=p^{-1}(p(x))$ is a smooth elliptic curve, and there is an unique $y\in \Sigma_x$ such that
$$ mx=y+k L_{\mid \Sigma_x}\,\,{\rm in}\,\,{\rm Pic}\,\Sigma_x.$$
This determines a rational map $\psi:\Sigma\dashrightarrow\Sigma,\,x\mapsto y$ which is of degree $m^2$.
Let $\tau:\widetilde{\Sigma}\rightarrow \Sigma$ be a birational morphism
such that $\psi\circ \tau$ is a morphism, and let
$$f:=\phi\circ\tau:\widetilde{\Sigma}\rightarrow S,\,g:=\phi\circ\psi\circ\tau:\widetilde{\Sigma}\rightarrow S.$$
\begin{rema}\label{remadujour}
{\rm The degree of $f$ is equal to the degree of $\phi$, hence independent of $m$.}
\end{rema}
\begin{lemm} \label{kcorrespk3} The image $\Sigma_m:=(f,g)(\widetilde{\Sigma})$ is a $K$-autocorrespondence of $S$,
which satisfies the following numerical properties:

1) For any $\eta\in H^{2,0}(S)$, $g^*\eta=mf^*\eta$.

2) $f_*g^* L=\lambda_m L\,{\rm in}\,\, {\rm Pic}\,S$, where
$$m\lambda_m\not\in\{0,\,m^2{\rm deg}\,f,\,m\,{\rm deg}\,f,\,{\rm deg}\,f\}$$ for $m$ large enough.

\end{lemm}

{\bf Proof.} By construction, we have for $\sigma\in \Sigma$
\begin{eqnarray}g(\sigma)=mf(\sigma)-kL^2\,\,{\rm in}\,\,CH_0(S).\label{numero4jan}
\end{eqnarray}
Thus $\Sigma_m$ is a $K$-correspondence and 1) is satisfied by Proposition \ref{criKcor}.

As $({\rm Pic}\,S)\otimes\mathbb{}=\mathbb{Q}L$, we certainly have a formula $f_*g^* L=\lambda_m L\,{\rm in}\,\, {\rm Pic}\,S$ and
it only remains to show that $m\lambda_m\not\in\{0,\,{\rm deg}\,f,\,m\,{{\rm deg}\,f},\,m^2{\rm deg}\,f\}$ for $m$ large. This is
however obvious,  as the degree of $f$ is independent of $m$ according to Remark \ref{remadujour}, while the intersection number
$g^*L\cdot\Sigma_b$ is equal to $2m^2d$, which implies that
the intersection number
$f_*\Sigma_b\cdot f_*g^*L=L\cdot f_*g^*L$ is $\geq 2m^2d$, so that $\lambda_m\geq m^2$.

\cqfd
 \begin{coro} \label{coro} For a very general pair $(S,L)$ as above, we have
 $$mf^*=g^*: H^2(S,\mathbb{Q})^{\perp c_1(L)}\rightarrow H^2(\widetilde{\Sigma},\mathbb{Q}).$$
\end{coro}
 {\bf Proof.} Indeed the morphism of Hodge structures $mf^*-g^*:H^2(S,\mathbb{Q})^{\perp c_1(L)}\rightarrow H^2(\widetilde{\Sigma},\mathbb{Q})$ vanishes on $H^{2,0}(S)$ by Lemma \ref{kcorrespk3}.
Its kernel $K$ is thus a Hodge substructure of $H^2(S,\mathbb{Q})^{\perp c_1(L)}$ which contains both $H^{2,0}(S)$ and
its complex conjugate $H^{0,2}(S)$. The orthogonal complement of $K$ in $H^2(S,\mathbb{Q})^{\perp c_1(L)}$
is thus contained in  ${\rm NS}\,(S)\otimes \mathbb{Q}$ and orthogonal to $c_1(L)$, hence is $0$
because  for a very general pair $(S,L)$, we have ${\rm NS}\,(S)\otimes \mathbb{Q}=\mathbb{Q}c_1(L)$.
\cqfd
\begin{coro}\label{coronew29dec}
 The eigenvalues of $f_*g^*$ acting on $H^*(S,\mathbb{Q})$ are
$${\rm deg }\, f,\,m\,{\rm deg } \,f,\,\lambda_m,\,{m^2}\,{\rm deg }\, f.$$
\end{coro}
{\bf Proof.} Indeed, $f_*g^*$ acts as ${\rm deg }\, f\,Id$ on $H^0(S,\mathbb{Q})$. Corollary
\ref{coro} and Proposition \ref{kcorrespk3}, 2) show that the eigenvalues of
$f_*g^*$ on $H^2(S,\mathbb{Q})$ are $m\,{\rm deg }\, f$ and $\lambda_m$, and finally
$f_*g^*$ acts as ${\rm deg }\, g\,Id$ on $H^4(S,\mathbb{Q})$. But ${\rm deg }\, g={m^2}{{\rm deg }\, f}$ because for  any non zero holomorphic $2$-form $\eta$ on $S$, we have
$g^*\eta=mf^*\eta$ and thus
$$\int_{\widetilde{\Sigma}} g^*\eta\wedge g^*\overline{\eta}={\rm deg }\, g\int_S\eta\wedge \overline{\eta}
=m^2\int_{\widetilde{\Sigma}} f^*\eta\wedge f^*\overline{\eta}=m^2{\rm deg }\, f\int_S\eta\wedge \overline{\eta},$$
where the integral $\int_S\eta\wedge \overline{\eta}$ is non zero.
\cqfd
We are  going to use now the above
constructions to  prove Theorem \ref{main14juilletintro}, (i) for families of $K3$ surfaces with generic Picard number $1$.

\vspace{0.5cm}

{\bf Proof of Theorem \ref{main14juilletintro}, (i).} We start with our family $\pi:\mathcal{S}\rightarrow B$ of $K3$ surfaces, which has the property that the very general fibers $\mathcal{S}_b$ have Picard number $1$. Let $\mathcal{L}$ be a relatively ample line bundle
on $\mathcal{S}$ of self-intersection $2d$. The construction mentioned previously of a $1$-dimensional family of irreducible elliptic curves with smooth total space
works in family, at least over a Zariski open set  of $B$.
Hence, replacing  $B$ by a Zariski open set, $\mathcal{S}$ by its inverse image under $\pi$, and $\mathcal{L}$ by a
rational  multiple
of $\mathcal{L}$
 if necessary, we can assume that there are a
family of smooth surfaces
$p:{\mathcal{T}}\rightarrow B$ and two morphisms
\begin{eqnarray}\label{newref}f,\,g: {\mathcal{T}}\rightarrow \mathcal{S}
\end{eqnarray}
whose fibers over $b\in B$ satisfy the conclusions of Lemma \ref{kcorrespk3} and Corollary \ref{coro}.

The relative cycle
\begin{eqnarray}
\label{numero}
\Gamma:=(f,g)_*({\mathcal{T}})+(\frac{m\,{\rm deg}\,f-\lambda_m}{2d})pr_1^*c_{1}(\mathcal{L})\cdot pr_2^*c_{1}(\mathcal{L})\in CH^2(\mathcal{S}\times_{B}\mathcal{S})_\mathbb{Q}
\end{eqnarray}
induces a morphism
$${\Gamma}^*:R\pi_{*}\mathbb{Q}\rightarrow R\pi_{*}\mathbb{Q},$$
which acts by Corollary \ref{coronew29dec} with respective eigenvalues
$$\lambda_0={\rm deg}\,f,\,\lambda_2=m{\rm deg}\,f,\, \lambda_4=m^2  {\rm deg}\,f$$ on $R^0\pi_{*}\mathbb{Q},\,R^2\pi_{*}\mathbb{Q},\,R^4\pi_{*}\mathbb{Q}$.

These three eigenvalues being distinct, we can apply
Lemma \ref{ledecomp} to the morphism $\Gamma^*$ acting on the object $R\pi_{*}\mathbb{Q}$
of the bounded derived category of sheaves of $\mathbb{Q}$-vector spaces on $B$.
We thus get a decomposition
\begin{eqnarray}\label{decompk323dec}
R\pi_{*}\mathbb{Q}=R^0\pi_{*}\mathbb{Q}\oplus R^2\pi_{*}\mathbb{Q}[-2]\oplus R^4\pi_{*}\mathbb{Q}[-4],
\end{eqnarray}
which is preserved by
$\Gamma^*$.
Note furthermore that $R^2\pi_{*}\mathbb{Q}[-2]$ is canonically the direct sum $\mathbb{Q}L[-2]\oplus R^2\pi_{*}\mathbb{Q}^{\perp L}[-2]$,
which provides us with the two direct summands
\begin{eqnarray}\label{summanddujour}\mathbb{Q}L[-2],\, R^2\pi_{*}\mathbb{Q}^{\perp L}[-2]
\end{eqnarray}
of $R\pi_{*}\mathbb{Q}$.

The proof of Theorem \ref{main14juilletintro} then  concludes with  the following:
\begin{prop}\label{prop23dec}  The decomposition (\ref{decompk323dec}) is multiplicative on a nonempty Zariski open set
of $B$.
\end{prop}
\cqfd

It remains to prove Proposition \ref{prop23dec}.
The proof will use the following lemma:
Let $f,\,g: \Sigma\rightarrow S$ be two morphisms from a smooth surface
$\Sigma$ to a $K3$ surface $S$ equipped with a line bundle $L$ with non zero self-intersection.
\begin{lemm}\label{lebeauvoi} Assume that
for some integers $m_1,\,m_2$, and for some fixed $0$-cycle $z_0$ of $S$, the relation
\begin{eqnarray}\label{eqn29dec}m_1f(\sigma)+m_2g(\sigma)=z_0
\end{eqnarray}
holds in $CH_0(S)$ for every $\sigma\in \Sigma$.
Then we have
\begin{eqnarray}\label{lechow} f_*g^*(c_1(L)^2)={\rm deg}\,g\, c_1(L)^2
\end{eqnarray}
in $CH_0(S)$.
\end{lemm}
{\bf Proof.} We just have to show that $ f_*g^*(c_1(L)^2)$ is proportional to $ c_1(L)^2$
in $CH_0(S)$, since $CH_0(S)$ has no torsion and the degrees of both sides in
(\ref{lechow})  are equal.
There are various criteria for a point $x$ of $S$ to be  proportional to
$c_1(L)^2$ in $CH_0(S)$. The one used in \cite{beauvoi} is that it is enough  that $x$ belongs
to some (singular) rational curve in $S$. The following criterion is a weaker characterization:

\begin{sublemm}\label{sublemm} Let $S$ be a $K3$ surface and $L$ be a line bundle on $S$
such that ${\rm deg}\,c_1(L)^2\not=0$.
 Let $j:C\rightarrow S$ be a non constant morphism from an irreducible curve $C$ to $S$, such that
$j_*:CH_0(C)\rightarrow CH_0(S)$ has for image $\mathbb{Z}$ (that is all points
$j(c),\,c\in C$, are rationally equivalent in $S$). Then for any $c\in C$,
$j(c)$ is  proportional to
$c_1(L)^2$ in $CH_0(S)$.
\end{sublemm}
{\bf Proof.} Let $H$ be an ample line bundle on $S$.
 As all points $j(c),\,c\in C$ are rationally equivalent in $S$, they are proportional in $CH_0(S)$
 to the cycle $j_*j^*H=j_*C\cdot H$, because the latter has a non zero degree. But it follows from
 Theorem \ref{BV}  that $j_*C\cdot H$ and $c_1(L)^2$ are proportional in $CH_0(S)$.
\cqfd
Coming back to our situation, we start from a singular rational curve $D\subset S$ in some ample linear system
$\mid H\mid$. Then we know by \cite[Thm 1]{beauvoi} that any point $x$ of $D$ is proportional
to $c_1(L)^2$ in $CH_0(S)$. On the other hand, the curve $g^{-1}(D)$ is connected and
 $f(g^{-1}(D))$ is not reduced to a point, because $f_*g^*H\not=0$ in ${\rm NS}(S)$. Let
 $C$ be a component of $g^{-1}(D)$ which is not contracted to a point by $f$. We now apply
 Sublemma  \ref{sublemm} to the morphism $f$ restricted to $C$. Indeed, as $g_*(c)$ is constant in $CH_0(S)$ because
 $g(C)$ is rational, it follows
 from  (\ref{eqn29dec}) that $f_*(c)$ is also constant in $CH_0(S)$.
 Hence $f_*(c)$ is proportional to $c_1(L)^2$ in $CH_0(S)$ by Sublemma \ref{sublemm}. As $g^{-1}(D)$ is connected, the same conclusion also holds
 for the components $C$ of $g^{-1}(D)$ which are contracted by $f$.
 As this is true for any $c\in C$, we get a fortiori that denoting by $g_C$ the restriction  of $g$ to $C$,
 $f_*g_C^*x$ is proportional to $c_1(L)^2$ in $CH_0(S)$.  Summing over all components $C$ of $g^{-1}(D)$, and recalling
  that $x$ is proportional to $c_1(L)^2$ in $CH_0(S)$ concludes the proof of  Lemma \ref{lebeauvoi}.
\cqfd
\begin{coro}\label{coro29dec} Over a nonempty Zariski open set of $B$, we have
\begin{eqnarray}\label{egalite29dec}\Gamma^*(c_1^{top}(\mathcal{L})^2)=m^2{\rm deg}\,f\,c_1^{top}(\mathcal{L})^2,
\end{eqnarray}
where $\Gamma$ is as in (\ref{numero}).

The  morphism
$c_1^{top}(\mathcal{L})^2\cup:\mathbb{Q}[-4]\rightarrow R\pi_*\mathbb{Q}$
factors through the direct summand $R^4\pi_*\mathbb{Q}[-4]$.
\end{coro}
{\bf Proof.} The second statement is an immediate consequence of
the first by definition of the decomposition.

Next, for any  point $b\in B$, $\Gamma_{b}^*$ acts as $f_*g^*$ on $CH_0(\mathcal{S}_b)$.
 Furthermore,  the pair $(f,g)$ satisfies the condition that
$$mf(\sigma)=g(\sigma)+kc_1(L)^2\,\,{\rm in}\,\,\, CH_0(\mathcal{S}_b)$$
for any $\sigma \in \mathcal{T}_b$.
As ${\rm deg}\,g=m^2{\rm deg}\,f$, Lemma \ref{lebeauvoi} tells us that
$$\Gamma_{b}^*(c_1(\mathcal{L}_b)^2)=f_*g^* (c_1(\mathcal{L}_b)^2)=m^2{\rm deg}\,f\,c_1(\mathcal{L}_b)^2$$
 in $CH_0(\mathcal{S}_b)$.

 The general principle \ref{genprinciple} then tells us that,
 for a nonempty Zariski open set $B^0$ of $B$,
 $$\Gamma^*(c_1^{top}(\mathcal{L})^2)=m^2{\rm deg}\,f\,c_1^{top}(\mathcal{L})^2\,\,{\rm in}\,\,H^4(\mathcal{S}^0, \mathbb{Q}).$$

\cqfd
\begin{coro}\label{remark29dec} The two morphisms  $\Gamma^*$ and  $f_*g^*$ agree,
over a nonempty Zariski open set of $B$, on the direct summand
$R^4\pi_{*}\mathbb{Q}[-4]$ of the decomposition (\ref{decompk323dec}). More precisely,  they both act by multiplication by
$m^2{\rm deg}\,f$ on this direct summand.

\end{coro}
{\bf Proof.} Indeed, this direct summand is equal by Corollary  \ref{coro29dec} to  the image of the morphism
\begin{eqnarray}\label{truc19jan}
\mathbb{Q}[-4]\rightarrow R\pi_{*}\mathbb{Q}
\end{eqnarray}
given by the class $c_1(\mathcal{L})^2$.
The difference
$f_*g^*-\Gamma^*$ is the morphism given by the class
$$\frac{m\,{\rm deg}\,f-\lambda_m}{2d}pr_1^*c_1^{top}(\mathcal{L})\cdot pr_2^*c_1^{top}(\mathcal{L})\in H^4(\mathcal{S}\times_B\mathcal{S},\mathbb{Q}),$$ hence is given up to a coefficient by the formula:
\begin{eqnarray}\label{formpresssee}pr_{1*}\circ ( pr_1^* c_1^{top}(\mathcal{L})\cup pr_2^*c_1^{top}(\mathcal{L})\cup)\circ pr_2^*: R\pi_{*}\mathbb{Q}\rightarrow R\pi_*\mathbb{Q}.
\end{eqnarray}
But the composition of the morphism (\ref{truc19jan})
with the morphism  (\ref{formpresssee})
obviously vanishes over a Zariski open set of $B$  because the class $c_1^{top}(\mathcal{L})^3\in H^6(\mathcal{S},\mathbb{Q})$ is the class of an algebraic cycle of codimension
  $3$.
\cqfd
We will also need the following easy lemma:
\begin{lemm}\label{easy21jan} 1) The morphisms $\Gamma_*$ and $f_*g^*$, restricted to the direct summand
$R^2\pi_{*}\mathbb{Q}^{\perp L}[-2]$ (see (\ref{summanddujour})), are equal.

2) The summand $\mathbb{Q}L[-2]$ of $R^2\pi_{*}\mathbb{Q}[-2]\subset R\pi_*\mathbb{Q}$ introduced in (\ref{summanddujour}) is
locally over $B$ in the Zariski topology generated by the class
$c_1^{top}(\mathcal{L})$, that is, is the image of the morphism
\begin{eqnarray}\label{cupL}c_1(\mathcal{L})\cup:\mathbb{Q}[-2]\rightarrow R\pi_*\mathbb{Q}.
\end{eqnarray}
\end{lemm}
{\bf Proof.} 1) Indeed their difference is  up to a coefficient the morphism  given by formula (\ref{formpresssee}).
But this
morphism
obviously
  vanishes on $ R^2\pi_{*}\mathbb{Q}^{\perp L}[-2]$, by the projection formula and because for degree reasons
it factors through the morphism of local systems
$$R^2\pi_{*}\mathbb{Q}
\stackrel{\cup c_1^{top}(\mathcal{L})}{\rightarrow}
R^4\pi_{*}\mathbb{Q}
\stackrel{\pi_{*}}{\rightarrow}\mathbb{Q}$$
which by definition vanishes on $R^2\pi_{*}\mathbb{Q}^{\perp L}$.

2)    Indeed, we  have
locally over $B$ in the Zariski topology
$$\Gamma^*c_1^{top}(\mathcal{L})=m\,{\rm deg}\,f\,c_1^{top}(\mathcal{L}).$$
  By definition of the decomposition, this implies  that locally over $B$, the morphism
  (\ref{cupL})
   takes value in the direct summand $R^2\pi_{*}\mathbb{Q}[-2]$ of  the decomposition. It then follows obviously that it locally belongs in fact to the direct summand
$\mathbb{Q}L[-2]$.
\cqfd
{\bf Proof of Proposition \ref{prop23dec}.}
We have the
data of the family of smooth surfaces
$p:\mathcal{T}\rightarrow B$ and of the morphisms
$f,\,g:\mathcal{T}\rightarrow \mathcal{S}$
as in
 (\ref{newref}).
 The  induced morphisms
 $$f^*:{R\pi}_*\mathbb{Q}\rightarrow Rp_*\mathbb{Q},\,g^*:{R\pi}_*\mathbb{Q}\rightarrow Rp_*\mathbb{Q},$$
  are multiplicative, i.e. compatible with cup-products on both sides.

Consider now our  decomposition
\begin{eqnarray}\label{rappeldecomp29dec}
 R\pi_{*}\mathbb{Q}\cong \oplus_i R^i\pi_{*}\mathbb{Q}[-i],
 \end{eqnarray}
together with the orthogonal decomposition of the local system $R^2\pi_{*}\mathbb{Q}$
$$R^2\pi_{*}\mathbb{Q}=R^2\pi_{*}\mathbb{Q}^{\perp L}\oplus \mathbb{Q}L.$$
The decomposition (\ref{rappeldecomp29dec}) is by definition
preserved by $\Gamma_{}^*$ and $\Gamma_{}^*$ acts with eigenvalues
$${\rm deg}\,f,\,m\,{\rm deg}\,f,\, m^2{\rm deg}\,f$$ on
the respective summands.

\begin{lemm}\label{lefg29dec} 1) Over a nonempty  Zariski open set of $B$, we have the equality
\begin{eqnarray}\label{truc19am}
g^*=mf^*: R^2\pi_{*}\mathbb{Q}^{\perp L}[-2]\rightarrow Rp_*\mathbb{Q}.
\end{eqnarray}
2) The morphism $f_*g^*: R\pi_{*}\mathbb{Q}\rightarrow R\pi_{*}\mathbb{Q}$ preserves the direct summand
$R^2\pi_{*}\mathbb{Q}^{\perp L}$ and acts by multiplication by $m{\rm deg}\,f$ on it.
\end{lemm}
{\bf Proof.}  2) follows from 1) by applying $f_*$ to both sides of (\ref{truc19am}).

To prove 1), note that the morphisms $f^*,\,g^*$ are induced by the classes of the
codimension $2$ cycles
$\Gamma_f:={\rm Graph}\,f,\,\Gamma_g:={\rm Graph}\,g$ in $\mathcal{T}\times_{B}\mathcal{S}_{}$.
For any $b\in B$, consider the cycle
$$\Gamma_b:=m\Gamma_{f,b}-\Gamma_{g,b}-k\,pr_2^*c_1(\mathcal{L}_b)^2\in CH^2(\mathcal{T}_b\times \mathcal{S}_b).$$
By construction, the induced map $\Gamma_{b*}: CH_0(\mathcal{T}_b)\rightarrow CH_0(\mathcal{S}_b)$ is equal to $0$.
It follows by applying the general principle \ref{genprinciple}
that, after passing to rational coefficients and modulo rational equivalence,
$\Gamma_b$ is supported on $D_b\times \mathcal{S}_b$ for some curve $D_b\subset \mathcal{T}_b$. However, as
${\rm Pic}\,^0(\mathcal{S}_b)=0$, denoting $\widetilde{D}_b$ the desingularization of ${D}_b$, we have
${\rm Pic}\,(\widetilde{D}_b\times \mathcal{S}_b)={\rm Pic}\,\widetilde{D}_b\oplus {\rm Pic}\,\mathcal{S}_b$.
We thus conclude that
\begin{eqnarray}\label{eqn8janv}
m\Gamma_{f,b}-\Gamma_{g,b}-kpr_2^*c_1(\mathcal{L}_b)^2=pr_1^*Z_b+ pr_1^*Z'_b\cdot pr_2^*Z''_b\,\,{\rm in}\,\,CH^2(\mathcal{T}_b\times \mathcal{S}_b)_{\mathbb{Q}},
\end{eqnarray}
for some zero cycle $Z_b\in CH^2(\mathcal{T}_b)$ and $1$-cycles $Z'_b$ on $\mathcal{T}_b$, $Z''_b$ on $\mathcal{S}_b$.
Note that the cycle $Z''_b$ has to be proportional to
$c_1(\mathcal{L}_b)$, since the point $b$ is general in $B$.

Applying again the general principle
\ref{genprinciple}, the pointwise equality (\ref{eqn8janv}) in the Chow groups of the fibers
produces the following  equality of cohomology classes over a Zariski open subset $B^0$:
\begin{eqnarray}\label{secondeqn8janv}m[\Gamma_{f}]-[\Gamma_{g}]-kpr_2^*c_1^{top}(\mathcal{L})^2
=pr_1^*[\mathcal{Z}]+ pr_1^*[\mathcal{Z}']\cup pr_2^*[\mathcal{Z}'']\,\,{\rm in}\,\,H^4(\mathcal{T}^0\times_{B^0} \mathcal{S}^0,{\mathbb{Q}})
\end{eqnarray}
for some codimension $2$ cycles $\mathcal{Z}\in CH^2(\mathcal{T}^0)_\mathbb{Q}$, and codimension $1$ cycles
$\mathcal{Z}'\in CH^1(\mathcal{T}^0)_\mathbb{Q},\,\mathcal{Z}''\in CH^1(\mathcal{S}^0)_\mathbb{Q}$, where  we may assume furthermore $[\mathcal{Z}'']=c_1^{top}(\mathcal{L})$ by shrinking $B^0$ if necessary.
We thus get  over $B^0$  the following equality of associated morphisms:
\begin{eqnarray}\label{thirdeqn8janv}mf^*=g^*+k(pr_2^*c_1^{top}(\mathcal{L})^2)^*+(pr_1^*[\mathcal{Z}]^*+ pr_1^*[\mathcal{Z}']\cup pr_2^*c_1^{top}(\mathcal{L}))^*:R\pi_{*}\mathbb{Q}\rightarrow Rp_*\mathbb{Q}.
\end{eqnarray}
The morphism
$$(pr_1^*[\mathcal{Z}])^*=pr_{1*}\circ ( pr_1^*[\mathcal{Z}]\cup)\circ pr_2^*: R\pi_{*}\mathbb{Q}\rightarrow Rp_*\mathbb{Q}$$
induced by the cycle class $pr_1^*[\mathcal{Z}]$ vanishes on $R^0\pi_{*}\mathbb{Q}\oplus R^2\pi_{*}\mathbb{Q}[-2]$, by the projection formula and because for degree reasons
$$pr_{1*}\circ pr_{2}^* =0,\,R^0\pi_{*}\mathbb{Q}\oplus R^2\pi_{*}\mathbb{Q}[-2]\rightarrow Rp_*\mathbb{Q}[-4]$$
vanishes.

Similarly, the
morphism
$$pr_{1*}\circ (  pr_1^*[\mathcal{Z}']\cup pr_2^*c_1^{top}(\mathcal{L})\cup)\circ pr_2^*: R\pi_{*}\mathbb{Q}\rightarrow Rp_*\mathbb{Q}$$
  vanishes on $R^0\pi_{*}\mathbb{Q}\oplus R^2\pi_{*}\mathbb{Q}^{\perp L}[-2]$, by the projection formula and because for degree reasons
it factors through the composite morphism
$$R^2\pi_{*}\mathbb{Q}
\stackrel{\cup c_1^{top}(\mathcal{L})}{\rightarrow}
R^4\pi_{*}\mathbb{Q}
\stackrel{\pi_{*}}{\rightarrow}\mathbb{Q}$$
which by definition vanishes on $R^2\pi_{*}\mathbb{Q}^{\perp L}$.

Using (\ref{thirdeqn8janv}), it only remains  to prove that
the restriction to $R^2\pi_{*}\mathbb{Q}^{\perp L}[-2]$
of the morphism induced by the class $pr_2^*c_1^{top}(\mathcal{L})^2$
$$pr_{1*}\circ(pr_2^*c_1^{top}(\mathcal{L})^2\cup)\circ pr_{2}^*:R\pi_{*}\mathbb{Q}\rightarrow Rp_*\mathbb{Q}$$  vanishes
over a Zariski open set of $B$.
Using (\ref{thirdeqn8janv}) and the above arguments, we conclude that on the direct summand
$R^2\pi_{*}\mathbb{Q}^{\perp L}[-2]$ and over a nonempty Zariski open set of $B$ we have
\begin{eqnarray}\label{thirdeqn21janv}mf^*=g^*+k(pr_2^*c_1^{top}(\mathcal{L})^2)^*:R^2\pi_{*}\mathbb{Q}^{\perp L}[-2]\rightarrow Rp_*\mathbb{Q}.
\end{eqnarray}
Applying $f_*$ to both sides, we conclude that
\begin{eqnarray}\label{thirdeqn41janv}m{\rm deg}\,f Id=f_*g^*+kf_*(pr_2^*c_1^{top}(\mathcal{L})^2)^*:R^2\pi_{*}\mathbb{Q}^{\perp L}[-2]\rightarrow R\pi_*\mathbb{Q}.
\end{eqnarray}
But $f_*g^*$  acts as $\Gamma_*$ on the direct summand $R^2\pi_{*}\mathbb{Q}^{\perp L}[-2]$
by Lemma \ref{easy21jan}, and by definition of the direct summand $R^2\pi_{*}\mathbb{Q}[-2]$, $\Gamma_*$ acts as
$m{\rm deg}\,f Id$ on it. Hence we have
$$f_*g^*=m{\rm deg}\,f Id :R^2\pi_{*}\mathbb{Q}^{\perp L}[-2]\rightarrow R\pi_*\mathbb{Q},$$
and comparing with (\ref{thirdeqn41janv}), we get that
\begin{eqnarray}\label{thirdeqn51janv} f_*\circ pr_2^*c_1^{top}(\mathcal{L})^2)^*=0:R^2\pi_{*}\mathbb{Q}^{\perp L}[-2]\rightarrow R\pi_*\mathbb{Q}.
\end{eqnarray}
It is now easy to see that the last equation implies
$$(pr_2^*c_1^{top}(\mathcal{L})^2)^*=0:R^2\pi_{*}\mathbb{Q}^{\perp L}[-2]\rightarrow Rp_*\mathbb{Q}.$$
Indeed, the morphism $(pr_2^*c_1^{top}(\mathcal{L})^2)^*:R\pi_*\mathbb{Q}\rightarrow Rp_*\mathbb{Q}$ factors as $p^*\circ \psi:\mathbb{Q}\rightarrow Rp_*\mathbb{Q}$, where
$\psi:R\pi_*\mathbb{Q}\rightarrow \mathbb{Q}$ is the composite morphism
$$R\pi_*\mathbb{Q}\stackrel{c_1^{top}(\mathcal{L})^2\cup}{\rightarrow } R\pi_*\mathbb{Q}[4]
\stackrel{\pi_*}{\rightarrow }\mathbb{Q},$$
and we have $f_*\circ p^*={\rm deg}\,f\circ \pi^*:\mathbb{Q}\rightarrow R\pi_*\mathbb{Q}$.
\cqfd
We now conclude the proof of Proposition \ref{prop23dec}.
Using Lemma \ref{lefg29dec}, we deduce now that,
in the decomposition (\ref{rappeldecomp29dec}), the cup-product map
$$\mu: R^2\pi_{*}\mathbb{Q}^{\perp L}[-2]\otimes R^2\pi_{*}\mathbb{Q}^{\perp L}[-2]
\rightarrow R\pi_{*}\mathbb{Q}$$
takes value in the direct summand $R^4\pi_{*}\mathbb{Q}[-4]$.
Indeed, we have
$g^*=mf^*$ on $ R^2\pi_{*}\mathbb{Q}^{\perp L}[-2]$ and thus
$$g^*\circ \mu=m^2f^*\circ \mu: R^2\pi_{*}\mathbb{Q}^{\perp L}[-2]\otimes R^2\pi_{*}\mathbb{Q}^{\perp L}[-2]\rightarrow Rp_*\mathbb{Q}.$$
Applying
$f_*$ on one hand, and taking the cup-product with $g^*c_1^{top}(\mathcal{L})$ on the other hand,  we conclude that, on
$R^2\pi_{*}\mathbb{Q}^{\perp L}[-2]\otimes R^2\pi_{*}\mathbb{Q}^{\perp L}[-2]$ we have:
\begin{eqnarray}\label{autretruc19jan}f_*g^*\circ \mu={\rm deg}\,f\,m^2 \mu:R^2\pi_{*}\mathbb{Q}^{\perp L}[-2]\otimes R^2\pi_{*}\mathbb{Q}^{\perp L}[-2]\rightarrow R\pi_{*}\mathbb{Q},\\
\label{autreautretruc}
g^*\circ \mu\circ( c_1^{top}(\mathcal{L})\cup)=m^2( g^*c_1^{top}(\mathcal{L})\cup)\circ f^*\circ \mu:\\
\nonumber R^2\pi_{*}\mathbb{Q}^{\perp L}[-2]\otimes R^2\pi_{*}\mathbb{Q}^{\perp L}[-2]\rightarrow Rp_{*}\mathbb{Q}[2],
\end{eqnarray}
Hence, by applying $f_*$ to the second equation (\ref{autreautretruc}), we get:
\begin{eqnarray}\label{latrois}f_*g^*\circ \mu\circ( c_1^{top}(\mathcal{L})\cup)=m^2\lambda_m \mu\circ( c_1^{top}(\mathcal{L})\cup):\\
\nonumber R^2\pi_{*}\mathbb{Q}^{\perp L}[-2]\otimes R^2\pi_{*}\mathbb{Q}^{\perp L}[-2]\rightarrow R\pi_{*}\mathbb{Q}[2].
\end{eqnarray}
Using Corollary \ref{remark29dec}, and Lemma \ref{lefg29dec}, 2), we get that
$f_*g^*$ preserves the decomposition (\ref{rappeldecomp29dec}), acting with
eigenvalues ${\rm deg}\,f$ on the first summand, $m{\rm deg}\,f$ and $\lambda_m$ on the summand
$R^2\pi_{*}[-2]$, and $m^2{\rm deg}\,f$ on the summand
$R^4\pi_{*}[-4]$. As $m^2\lambda_m\not\in\{{\rm deg}\,f ,\,m{\rm deg}\,f,\,\lambda_m,\,m^2{\rm deg}\,f\}$
 by Lemma
\ref{kcorrespk3}, 2), we first conclude from (\ref{latrois}) that
$$ \mu\circ( c_1^{top}(\mathcal{L})\cup)=c_1^{top}(\mathcal{L})\cup \circ \mu$$
vanishes on $R^2\pi_{*}\mathbb{Q}^{\perp L}[-2]\otimes R^2\pi_{*}\mathbb{Q}^{\perp L}[-2]$.

 Next,
we conclude from (\ref{autretruc19jan}) that $\mu:R^2\pi_{*}\mathbb{Q}^{\perp L}[-2]\otimes R^2\pi_{*}\mathbb{Q}^{\perp L}[-2]\rightarrow R\pi_{*}\mathbb{Q}$ takes value in the direct summand with is the sum $\mathbb{Q}L[-2]\oplus R^4\pi_{*}\mathbb{Q}[-4]$
(the second summand being possible if $\lambda_m=m^2{\rm deg}\,f$).
However, as its composition  with the cup-product map $ c_1^{top}(\mathcal{L})\cup$ vanishes, we easily conclude that
it actually takes value in the summand $ R^4\pi_{*}\mathbb{Q}[-4]$, because the cup-product map $ c_1^{top}(\mathcal{L})\cup$ induces an isomorphism
$\mathbb{Q}L[-2]\cong R^4\pi_{*}\mathbb{Q}[-2]$, as follows from Lemma \ref{easy21jan}, 2) and Corollary
\ref{coro29dec}.

It remains to see what happens on the other summands :
First of all, Lemma \ref{easy21jan}, 2) says that the summand $\mathbb{Q}L[-2]$ of $R^2\pi_{*}\mathbb{Q}[-2]$  is,
 over a nonempty Zariski open subset $B$, the image of the morphism
$c_1^{top}(\mathcal{L})\cup:\mathbb{Q}[-2]\rightarrow R\pi_*\mathbb{Q}$. On the other hand,
Corollary  \ref{coro29dec} says that   the direct summand $R^4\pi_{*}\mathbb{Q}[-4]$ is over a nonempty Zariski
open set $B^0$ of $B$ the image of the morphism
$$c_1^{top}(\mathcal{L})^2:\mathbb{Q}[-4]\rightarrow R\pi_*\mathbb{Q}.$$
It follows immediately  that
 for the summand $\mathbb{Q}L[-2]={\rm Im}\,c_1^{top}(\mathcal{L})\cup$
 $$ \mu :\mathbb{Q}L[-2]\otimes \mathbb{Q}L[-2]\rightarrow  R\pi_{*}\mathbb{Q}$$
 takes value on $B^0$ in the direct summand $R^4\pi_{*}\mathbb{Q}[-4]$.

 Consider now the cup-product
 $$ R^2\pi_{*}\mathbb{Q}^{\perp L}\otimes \mathbb{Q}L[-2]\rightarrow R\pi_{*}\mathbb{Q}.$$
 We claim that it vanishes over a nonempty Zariski open set of $B$.

  Indeed, Lemma \ref{lefg29dec}  tells that over a nonempty Zariski open set of $B$,
 $$g^*=mf^*: R^2\pi_{*}\mathbb{Q}^{\perp L}[-2]\rightarrow Rp_*\mathbb{Q}.$$

 It follows that
 $$ g^*\circ\mu=\mu \circ (g^*\otimes g^*) =\mu \circ (mf^*\otimes g^*) :R^2\pi_{*}\mathbb{Q}^{\perp L}[-2]\otimes \mathbb{Q}L[-2]\rightarrow Rp_*\mathbb{Q}.$$
 Applying  the projection formula, we get  that
 $$f_*g^*\circ\mu=\mu\circ (m Id\otimes f_*g^*):R^2\pi_{*}\mathbb{Q}^{\perp L}[-2]\otimes \mathbb{Q}L[-2]\rightarrow R\pi_{*}\mathbb{Q}.$$

On the other hand, we know by Lemma \ref{easy21jan}, 2) that $f_*g^*$ sends, locally over $B$, the summand
  $\mathbb{Q}L[-2]$ to itself, acting on it by multiplication by $\lambda_m$.
  It follows that
  $$f_*g^*\circ\mu=m\lambda_m\mu :R^2\pi_{*}\mathbb{Q}^{\perp L}[-2]\otimes \mathbb{Q}L[-2]\rightarrow R\pi_{*}\mathbb{Q},$$
  and finally we conclude that $f_*g^*\circ\mu=0$ on $R^2\pi_{*}\mathbb{Q}^{\perp L}[-2]\otimes \mathbb{Q}L[-2]$ because
  $m\lambda_m$ is not an eigenvalue of $f_*g^*$ acting on the cohomology of $R\pi_{*}\mathbb{Q}$
  by Corollary \ref{coronew29dec} and  Lemma
\ref{kcorrespk3}, 2).

  To conclude the proof of the multiplicativity, we just have to check that
  the cup-product map vanishes  over a nonempty Zariski open subset of $B$ on the
  summands $R^2\pi_{*}\mathbb{Q}[-2]\otimes  R^4\pi_{*}\mathbb{Q}[-4]$ and
  $R^4\pi_{*}\mathbb{Q}[-4]\otimes R^4\pi_{*}\mathbb{Q}[-4]$. The proof works exactly as before, by an eigenvalue computation for
  the summand $R^2\pi_{*}\mathbb{Q}^{\perp L}[-2]\otimes  R^4\pi_{*}[-4]$. For the other terms,
  this is clear because we have seen that over an adequate Zariski open subset of $B$,
  the factors  are  generated by  classes $c_1^{top}(\mathcal{L})$, $c_1^{top}(\mathcal{L})^2$,
  whose products are classes of algebraic cycles on $\mathcal{S}$
  of codimension at least $3$, hence vanishing  over a nonempty Zariski open subset of $B$.

\cqfd

\subsection{Alternative proof}
In this section we give a different proof of Theorem
\ref{main14juilletintro}, which also provides a proof of the second statement (ii). It heavily uses
the following  result proved in \cite[Proposition 3.2]{beauvoi}, whose proof is rather intricate.

   \begin{theo} (Beauville-Voisin 2004)\label{decompdiagonalpetiteintro}
Let $S$ be a smooth projective $K3$ surface, $L$ an ample line bundle on $S$ and $o:=\frac{1}{{\rm deg}_{S}\,c_1(L)^2}L^2\in CH^2(S)_\mathbb{Q}$. We have
\begin{eqnarray}\label{equadeltapetiteintro}\Delta=\Delta_{12}\cdot o_3+(perm.)-  (o_1\times o_2\times S+(perm.))\,\,{\rm in}\,\,CH^4(S\times S\times S)_\mathbb{Q}.
\end{eqnarray}
\end{theo}
(We recall that  ``$+(perm.)$'' means that we symmetrize the considered expression in the indices. The lower index $i$ means ``pull-back of the considered cycle under
the $i$-th projection $S^3\rightarrow S$'', and the lower index $ij$ means ``pull-back of the considered cycle under
the projection $S^3\rightarrow S^2$ onto the product of the $i$th and $j$-th factor''.

\vspace{0.5cm}

 {\bf Second proof of Theorem \ref{main14juilletintro}.}
Let us choose
 a relatively ample line bundle
 $L$ on $\mathcal{X}$, and let
 $$o_\mathcal{X}:=\frac{1}{{\rm deg}_{\mathcal{X}_t}\,c_1(L)^2}L^2\in CH^2(\mathcal{X})_\mathbb{Q}.$$
 By Theorem \ref{BV} and the general principle
 \ref{genprinciple}, this  cycle, which is of relative degree $1$, does not depend on the choice of
  $L$ up to shrinking the base $B$.
The cohomology classes
   $$pr_1^*[o_\mathcal{X}]=[Z_0],\,pr_2^*[o_\mathcal{X}]=[Z_4]\in H^4(\mathcal{X}\times_B \mathcal{X},\mathbb{Q})$$
   of the two codimension $2$ cycles $Z_0:=pr_1^*o_\mathcal{X}$ and $Z_4:=pr_2^*o_\mathcal{X}$, where
   $pr_i:\mathcal{X}\times_B \mathcal{X}\rightarrow B$
   are the two projections, provide
   morphisms in the derived category:
   $$P_0: R\pi_*\mathbb{Q}\rightarrow R\pi_*\mathbb{Q},\,\,
   P_4:R\pi_*\mathbb{Q}\rightarrow R\pi_*\mathbb{Q}$$
   \begin{eqnarray}
   \label{formP1*}P_0:=pr_{2*}\circ ([Z_0]\cup)\circ pr_1^*,\,P_4:=pr_{2*}\circ (pr_2^*[Z_4]\cup)\circ pr_1^*.
   \end{eqnarray}
   \begin{lemm} \label{le29jui2011}
   (i) The morphisms $P_0$, $P_4$ are projectors of $R\pi_*\mathbb{Q}$.

   (ii)  $P_0\circ P_4=P_4\circ P_0=0$ over a Zariski dense open set of $B$.
   \end{lemm}
   {\bf Proof.} (i) We compute $P_0\circ P_0$. From (\ref{formP1*}) and the projection formula
   \cite[Prop. 8.3]{fulton}, we get
   that $P_0\circ P_0$ is the morphism $R\pi_*\rightarrow R\pi_*$ induced
   by the following cycle class
   \begin{eqnarray}\label{classeidiote29juillet}
   p_{13*}(p_{12}^*[Z_0]\cup p_{23}^*[Z_0])\in H^4(\mathcal{X}\times_B \mathcal{X},\mathbb{Q}),
   \end{eqnarray}
   where the $p_{ij}$ are the various projections from
   $\mathcal{X}\times_B \mathcal{X}\times_B \mathcal{X}$ to $\mathcal{X}\times_B \mathcal{X}$.
   We use now the fact that
   $p_{12}^*[Z_0]=p_1^*[o_\mathcal{X}],\,\,p_{23}^*[Z_0]=p_2^*[o_\mathcal{X}]$,
    where the $p_{i}$'s are the various projections from
   $\mathcal{X}\times_B \mathcal{X}\times_B \mathcal{X}$ to $\mathcal{X}$, so that
   (\ref{classeidiote29juillet}) is equal to
   \begin{eqnarray}\label{classeidiote29juillet2}p_{13*}(p_1^*[o_\mathcal{X}]\cup p_2^*[o_\mathcal{X}]).
   \end{eqnarray}
   Using the projection formula, this class
   is equal
   to
   $$pr_1^*[o_\mathcal{X}]\cup pr_2^*(\pi_*[o_\mathcal{X}])=pr_1^*[o_\mathcal{X}]\cup pr_2^*(1_B)=pr_1^*[o_\mathcal{X}]=[Z_1].$$
    This completes the proof for $P_0$ and exactly the same proof works  for $P_4$.

   (ii) We compute $P_0\circ P_4$ : From (\ref{formP1*}) and  the projection formula
   \cite[Prop. 8.3]{fulton}, we get
   that $P_0\circ P_4$ is the morphism $R\pi_*\rightarrow R\pi_*$ induced
   by the following cycle class
   \begin{eqnarray}\label{autreclasseidiote29juillet}
   p_{13*}(p_{12}^*[Z_4]\cup p_{23}^*[Z_0])\in H^4(\mathcal{X}\times_B \mathcal{X},\mathbb{Q}),
   \end{eqnarray}
   where the $p_{ij}$ are the various projections from
   $\mathcal{X}\times_B \mathcal{X}\times_B \mathcal{X}$ to $\mathcal{X}\times_B \mathcal{X}$.
    We use now the fact that
   $p_{12}^*[Z_4]=p_2^*[o_\mathcal{X}],\,\,p_{23}^*[Z_0]=p_2^*[o_\mathcal{X}]$,
    where the $p_{i}$'s are the various projections from
   $\mathcal{X}\times_B \mathcal{X}\times_B \mathcal{X}$ to $\mathcal{X}$, so that
   (\ref{autreclasseidiote29juillet}) is equal to
   \begin{eqnarray}\label{autreclasseidiote29juillet2}p_{13*}(p_2^*[o_\mathcal{X}]\cup p_2^*[o_\mathcal{X}]).
   \end{eqnarray}
   But the class  $p_2^*[o_\mathcal{X}]\cup p_2^*[o_\mathcal{X}]=p_2^*([o_\mathcal{X}\cdot o_\mathcal{X}])$
   vanishes over a Zariski dense open set of $B$ since the cycle
   $o_\mathcal{X}\cdot o_\mathcal{X}$ has codimension $4$ in $\mathcal{X}$. This shows that
   $P_0\circ P_4=0$ over a Zariski dense open set of $B$
   and the proof for $P_4\circ P_0$ works in the same way.

   \cqfd
   Using Lemma \ref{le29jui2011}, we get (up to passing to
   a Zariski dense open set of $B$) a third projector
   $$P_2:=Id-P_0-P_4$$
   acting on $R\pi_*\mathbb{Q}$ and commuting with the two other ones.

   It is well-known (cf. \cite{murre}) that the action of these three projectors on cohomology
   are given by
   $$P_{0}=0\,\,{\rm on}\,\, R^2\pi_*\mathbb{Q},\,R^4\pi_*\mathbb{Q},\,\,\,\,P_{0*}=Id\,\,{\rm on}\,\, R^0\pi_*\mathbb{Q}.$$
   $$P_{4}=0\,\,{\rm on}\,\, R^2\pi_*\mathbb{Q},\,R^0\pi_*\mathbb{Q},\,\,\,\,P_{4*}=Id\,\,{\rm on}\,\, R^4\pi_*\mathbb{Q}.$$
   $$P_2=0\,\,{\rm on}\,\, R^0\pi_*\mathbb{Q},\,\,R^4\pi_*\mathbb{Q},\,\,\,\,P_{2}=Id\,\,{\rm on}\,\, R^2\pi_*\mathbb{Q}.$$

   As a consequence, we get (for example using Lemma \ref{ledecomp}) a decomposition
   \begin{eqnarray}\label{ladecomp29ju2011}
   R\pi_*\mathbb{Q}\cong \oplus R^i\pi_*\mathbb{Q}[-i],
   \end{eqnarray}
   where the corresponding projectors  of $R\pi_*\mathbb{Q}$
    identify respectively to $P_0,\,P_2,\,P_4$.

We now prove the following result,
\begin{prop} \label{aprouverlabarbe}
Assume the cohomology class of the relative small diagonal $\Delta\subset \mathcal{X}\times_B\mathcal{X}\times_B\mathcal{X}$
satisfies the equality
\begin{eqnarray} \label{equasmalldiagrel29ju}
[\Delta]=p_1^*[o_\mathcal{X}]\cup p_{23}^*[\Delta_\mathcal{X}]+(perm.)-(p_1^*[o_\mathcal{X}]\cup p_2^*[o_\mathcal{X}] +(perm.))
\end{eqnarray}
where the $p_{ij},\,p_i$'s are as above and $\Delta_\mathcal{X}$ is the relative
diagonal $\mathcal{X}\subset \mathcal{X}\times_B\mathcal{X}$,
then, over some Zariski dense open set $B^0\subset B$, we have:

(i) The decomposition (\ref{ladecomp29ju2011}) is multiplicative.

(ii)  The class of the diagonal $[\Delta_\mathcal{X}]\in H^4(\mathcal{X}\times_B \mathcal{X},\mathbb{Q})$
belongs to the direct summand $$H^0(B,R^4(\pi,\pi)_*\mathbb{Q})\subset H^4(\mathcal{X}\times_B \mathcal{X},\mathbb{Q})$$ induced by the decomposition
(\ref{ladecomp29ju2011}).
\end{prop}
Admitting Proposition \ref{aprouverlabarbe}, the end of the proof of Theorem \ref{main14juilletintro}
is as follows:
By Theorem \ref{decompdiagonalpetiteintro}, we know that the relation
$$\Delta_t=p_1^*o_{\mathcal{X}_t}\cdot p_{23}^*\Delta_{\mathcal{X}_t}+(perm.)-(p_1^*o_{\mathcal{X}_t}\cdot p_2^*o_{\mathcal{X}_t} +(perm.))$$
holds in $CH_2(\mathcal{X}_t\times \mathcal{X}_t\times \mathcal{X}_t,\mathbb{Q})$ for any $t\in B$.
By the general principle \ref{genprinciple}, we conclude that there exists a Zariski dense open set
$B^0$ of $B$ such that
(\ref{equasmalldiagrel29ju}) holds in $H^8(\mathcal{X}\times_B\mathcal{X}\times_B\mathcal{X},\mathbb{Q})$.
The statements (i) and (ii) of Theorem \ref{main14juilletintro} thus follow respectively from the statements (i) and (ii)  of Proposition \ref{aprouverlabarbe}. As proved in Lemma \ref{newlemma10aout},
the statement  (iii) of  Theorem
\ref{main14juilletintro} is implied by  (i).
   \cqfd
{\bf Proof of Proposition \ref{aprouverlabarbe}.}
(i) We want to show that
$$P_k\circ \cup\circ (P_i\otimes P_j):R\pi_*\mathbb{Q}\otimes R\pi_*\mathbb{Q}\rightarrow R\pi_*\mathbb{Q}$$
vanishes for $k\not=i+j$.

We note that
$$\cup:R\pi_*\mathbb{Q}\otimes R\pi_*\mathbb{Q}\rightarrow R\pi_*\mathbb{Q}$$
is induced, via the relative K\"{u}nneth decomposition $$R\pi_*\mathbb{Q}\otimes R\pi_*\mathbb{Q}\cong R(\pi,\pi)_*\mathbb{Q}$$
by the class $[\Delta]$ of the small relative diagonal
in $\mathcal{X}\times_B\mathcal{X}\times_B\mathcal{X}$, seen as a relative correspondence between
$\mathcal{X}\times_B\mathcal{X}$ and $\mathcal{X}$, while $P_0,\,P_4,\,P_2$ are induced by the cycle classes
$[Z_0],\,[Z_4],\,[Z_2]\in H^4(\mathcal{X}\times_B\mathcal{X},\mathbb{Q})$, where $Z_2:=\Delta_\mathcal{X}-Z_0-Z_4\subset \mathcal{X}\times_B\mathcal{X}$.
It thus suffices to show
that the cycle classes
$$[ Z_4\circ \Delta\circ (Z_0\times_B Z_0)],\,[ Z_2\circ \Delta\circ (Z_0\times_B Z_0)],$$
$$[ Z_0\circ \Delta\circ (Z_2\times_B Z_2)],\,[ Z_2\circ \Delta\circ (Z_2\times_B Z_2)],$$
$$[ Z_0\circ \Delta\circ (Z_4\times_B Z_4)],\,[ Z_2\circ \Delta\circ (Z_4\times_B Z_4)],\,[ Z_4\circ \Delta\circ (Z_4\times_B Z_4)],$$
$$[ Z_0\circ \Delta\circ (Z_0\times_B Z_4)],\,[ Z_2\circ \Delta\circ (Z_0\times_B Z_4)],$$
$$[ Z_0\circ \Delta\circ (Z_2\times_B Z_4)],\,[ Z_2\circ \Delta\circ (Z_2\times_B Z_4)],\,[ Z_4\circ \Delta\circ (Z_2\times_B Z_4)],$$
$$ [ Z_0\circ \Delta\circ (Z_0\times_B Z_2)],\,[ Z_4\circ \Delta\circ (Z_0\times_B Z_2)],$$
vanish in $H^8(\mathcal{X}\times_B\mathcal{X}\times_B\mathcal{X},\mathbb{Q})$
over a dense Zariski open set of $B$. Here, all the compositions of correspondences are  over $B$.
Equivalently, it suffices to prove the following equality of cycle classes
in $H^8(\mathcal{X}^0\times_B\mathcal{X}^0\times_B\mathcal{X}^0,\mathbb{Q})$, $\mathcal{X}^0=\pi^{-1}(B^0)$, for  a Zariski dense open set of $B^0$ of $B$:
\begin{eqnarray}\label{equanou3aout2011}
[\Delta]=[ Z_0\circ \Delta\circ (Z_0\times_B Z_0)]+[ Z_4\circ \Delta\circ (Z_2\times_B Z_2)]
+[ Z_2\circ \Delta\circ (Z_0\times_B Z_2)]\\
\nonumber
+[ Z_2\circ \Delta\circ (Z_2\times_B Z_0)]+[ Z_4\circ \Delta\circ (Z_0\times_B Z_4)]+[ Z_4\circ \Delta\circ (Z_4\times_B Z_0)].
\end{eqnarray}
Replacing $Z_2$ by $\Delta_\mathcal{X}-Z_0-Z_4$, we get
$$Z_2\times_B Z_2=\Delta_\mathcal{X}\times_B \Delta_\mathcal{X}-\Delta_\mathcal{X}\times_B Z_0-\Delta_\mathcal{X}\times_B Z_4-Z_0\times_B\Delta_\mathcal{X}$$
$$-Z_4\times_B \Delta_\mathcal{X}+Z_0\times_B Z_0+Z_4\times_B Z_4+Z_0\times_B Z_4+Z_4\times_B Z_0$$ and thus
(\ref{equanou3aout2011}) becomes
\begin{eqnarray}\label{equanou3aout2011aprem1}
[\Delta]=[ Z_0\circ \Delta\circ (Z_0\times_B Z_0)]+[ Z_4\circ \Delta\circ (\Delta_\mathcal{X}\times_B \Delta_X)]\,\,\,\\
\nonumber
-[ Z_4\circ \Delta\circ (\Delta_\mathcal{X}\times_B Z_0)]-[ Z_4\circ \Delta\circ (\Delta_\mathcal{X}\times_B Z_4)]
-[ Z_4\circ \Delta\circ (Z_0\times_B\Delta_\mathcal{X})]\,\,\,
\\
\nonumber
-[ Z_4\circ \Delta\circ (Z_4\times_B\Delta_\mathcal{X})]
+[ Z_4\circ \Delta\circ (Z_0\times_B Z_0)]+[ Z_4\circ \Delta\circ (Z_4\times_B Z_4)]\,\,\,
\\
\nonumber
+[ Z_4\circ \Delta\circ (Z_0\times_B Z_4)]+[ Z_4\circ \Delta\circ (Z_4\times_B Z_0)]
+[ Z_2\circ \Delta\circ (Z_0\times_B \Delta_\mathcal{X})]\,\,\,\\
\nonumber
-[ Z_2\circ \Delta\circ (Z_0\times_B Z_0)]-[ Z_2\circ \Delta\circ (Z_0\times_B Z_4)]
+[ Z_2\circ \Delta\circ (\Delta_\mathcal{X}\times_B Z_0)]\,\,\,
\\
\nonumber
-[ Z_2\circ \Delta\circ (Z_0\times_B Z_0)]
-[ Z_2\circ \Delta\circ (Z_4\times_B Z_0)]
+[ Z_4\circ \Delta\circ (Z_0\times_B Z_4)]\,\,\,
\\
\nonumber
+[ Z_4\circ \Delta\circ (Z_4\times_B Z_0)].\,\,\,
\end{eqnarray}
We now have the following lemma:
\begin{lemm}\label{dernierlemme3aout2011} We have the following equalities of cycles in
$CH^4(\mathcal{X}\times_B\mathcal{X}\times_B\mathcal{X})_\mathbb{Q}$ (or relative correspondences between $\mathcal{X}\times_B\mathcal{X}$ and $\mathcal{X}$)
\begin{eqnarray}\label{aout3eq1}\Delta\circ (Z_0\times_B Z_0)=p_1^*o_\mathcal{X}\cdot p_2^*o_\mathcal{X},\end{eqnarray}
\begin{eqnarray}\label{aout3eq2}\Delta\circ (\Delta_\mathcal{X}\times_B \Delta_\mathcal{X})=\Delta,\end{eqnarray}
\begin{eqnarray}\label{aout3eq3}\Delta\circ (\Delta_\mathcal{X}\times_B Z_0)=p_{13}^*\Delta_\mathcal{X}\cdot p_2^*o_\mathcal{X},\end{eqnarray}
\begin{eqnarray}\label{aout3eq4} \Delta\circ (\Delta_\mathcal{X}\times_B Z_4)=p_1^*o_\mathcal{X}\cdot p_3^*o_\mathcal{X},\end{eqnarray}
\begin{eqnarray}\label{aout3eq5} \Delta\circ (Z_0\times_B\Delta_\mathcal{X})=p_1^*o_\mathcal{X}\cdot p_{23}^*\Delta_\mathcal{X},\end{eqnarray}
\begin{eqnarray}\label{aout3eq6}\Delta\circ (Z_4\times_B\Delta_\mathcal{X})=p_{2}^*o_\mathcal{X}\cdot p_3^*o_\mathcal{X},\end{eqnarray}
\begin{eqnarray}\label{aout3eq7}\Delta\circ (Z_4\times_B Z_4)=p_3^*(o_\mathcal{X}\cdot o_\mathcal{X}),\end{eqnarray}
\begin{eqnarray}\label{aout3eq8}\Delta\circ (Z_0\times_B Z_4)=p_1^*o_\mathcal{X}\cdot p_3^*o_\mathcal{X},\end{eqnarray}
\begin{eqnarray}\label{aout3eq9} \Delta\circ (Z_4\times_B Z_0)=p_2^*o_\mathcal{X}\cdot p_3^*o_\mathcal{X},\end{eqnarray}
where the $p_i$'s, for $i=1,2,3$ are the projections from
$\mathcal{X}\times_B\mathcal{X}\times_B\mathcal{X}$ to $\mathcal{X}$ and the $p_{ij}$ are the projections from
$\mathcal{X}\times_B\mathcal{X}\times_B\mathcal{X}$ to $\mathcal{X}\times_B\mathcal{X}$.
\end{lemm}
{\bf Proof.} Equation (\ref{aout3eq2}) is obvious.
Equations (\ref{aout3eq1}), (\ref{aout3eq7}), (\ref{aout3eq8}), (\ref{aout3eq9}) are all similar. Let us just prove
(\ref{aout3eq8}). The cycle $Z_4$ is $\mathcal{X}\times_Bo_\mathcal{X}\subset \mathcal{X}\times_B\mathcal{X}$,
and similarly $Z_0=o_\mathcal{X}\times_B\mathcal{X}\subset \mathcal{X}\times_B\mathcal{X}$, hence
$Z_0\times_B Z_4$ is the cycle
\begin{eqnarray}\label{formule211aout}
\{(o_{\mathcal{X}_b},x,y,o_{\mathcal{X}_b}),\,x\in \mathcal{X}_b,\,y\in \mathcal{X}_b,\,b\in B\}\subset \mathcal{X}\times_B\mathcal{X}\times_B\mathcal{X}\times_B\mathcal{X}.
\end{eqnarray}
(It turns out that in this case, we do not have to take care about the ordering we take for the last inclusion.)
Composing over $B$ with $\Delta\subset \mathcal{X}\times_B\mathcal{X}\times_B\mathcal{X}$ is done by
 taking the pull-back of (\ref{formule211aout}) under
$p_{1234}:\mathcal{X}^{5/B}\rightarrow \mathcal{X}^{4/B}$,
intersecting with $p_{345}^*\Delta$, and projecting the resulting cycle to $\mathcal{X}^{3/B}$
via $p_{125}$. The resulting cycle is obviously
$$\{(o_{\mathcal{X}_b},x,o_{\mathcal{X}_b}),\,x\in \mathcal{X}_b,\,b\in B\}\subset \mathcal{X}\times_B\mathcal{X}\times_B\mathcal{X},$$
which proves (\ref{aout3eq8}).

For the last formulas which are all of the same kind, let us  just prove
(\ref{aout3eq3}).
Recall that $Z_0=o_\mathcal{X}\times_B\mathcal{X}\subset \mathcal{X}\times_B\mathcal{X}$.
Thus $\Delta_\mathcal{X}\times_B Z_0$ is the cycle
$$\{(x,x,o_{\mathcal{X}_b},y),\,x\in \mathcal{X}_b,\,y\in \mathcal{X}_b,\,b\in B\}\subset \mathcal{X}\times_B\mathcal{X}\times_B\mathcal{X}\times_B\mathcal{X}.$$
But we have to see this cycle as a relative self-correspondence of
$\mathcal{X}\times_B\mathcal{X}$,
for which the right ordering is
\begin{eqnarray}
\label{effroi3aout}\{(x,o_{\mathcal{X}_b},x,y),\,x\in \mathcal{X}_b,\,y\in \mathcal{X}_b,\,b\in B\}\subset \mathcal{X}\times_B\mathcal{X}\times_B\mathcal{X}\times_B\mathcal{X}.
\end{eqnarray}
Composing over $B$ with $\Delta\subset \mathcal{X}\times_B\mathcal{X}\times_B\mathcal{X}$ is done again by
 taking the pull-back of (\ref{effroi3aout}) by
$p_{1234}:\mathcal{X}^{5/B}\rightarrow \mathcal{X}^{4/B}$,
intersecting with $p_{345}^*\Delta$, and projecting the resulting cycle to $\mathcal{X}^{3/B}$
via $p_{125}$.
Since $\Delta=\{(z,z,z),\,z\in \mathcal{X}$, the considered intersection
is $\{(x,o_{\mathcal{X}_b},x,x,x),\,x\in \mathcal{X}_b,\,\,b\in B\}$, and thus the projection via
$p_{125}$ is $\{(x,o_{\mathcal{X}_b},x),\,x\in \mathcal{X}_b,\,\,b\in B\}$, thus proving
(\ref{aout3eq3}).

\cqfd
Using Lemma \ref{dernierlemme3aout2011} and the fact that
the cycle $p_3^*(o_\mathcal{X}\cdot o_\mathcal{X})$  vanishes by dimension reasons over
a dense Zariski open set of $B$, (\ref{equanou3aout2011aprem1}) becomes, after passing to
a Zariski open set of $B$ if necessary:

\begin{eqnarray}
[\Delta]=[ Z_0\circ (p_1^*o_\mathcal{X}\cdot p_2^*o_\mathcal{X})]+[ Z_4\circ \Delta]\,\,\,\\
\nonumber
-[ Z_4\circ (p_{13}^*\Delta_\mathcal{X}\cdot p_2^*o_\mathcal{X})]-[ Z_4\circ (p_1^*o_\mathcal{X}\cdot p_3^*o_\mathcal{X})]
-[ Z_4\circ (p_1^*o_\mathcal{X}\cdot p_{23}^*\Delta_\mathcal{X})]\,\,\,
\\
\nonumber
-[ Z_4\circ (p_{2}^*o_\mathcal{X}\cdot p_3^*o_\mathcal{X})]
+[ Z_4\circ (p_1^*o_\mathcal{X}\cdot p_2^*o_\mathcal{X})]
+[ Z_4\circ (p_1^*o_\mathcal{X}\cdot p_3^*o_\mathcal{X})]\,\,\,\\
\nonumber+[ Z_4\circ (p_2^*o_\mathcal{X}\cdot p_3^*o_\mathcal{X})]
+[ Z_2\circ (p_1^*o_\mathcal{X}\cdot p_{23}^*\Delta_\mathcal{X})]
-[ Z_2\circ (p_1^*o_\mathcal{X}\cdot p_2^*o_\mathcal{X})]\,\,\,
\\
\nonumber-[ Z_2\circ (p_1^*o_\mathcal{X}\cdot p_3^*o_\mathcal{X})]
+[ Z_2\circ (p_{13}^*\Delta_\mathcal{X}\cdot p_2^*o_\mathcal{X})]
-[ Z_2\circ (p_1^*o_\mathcal{X}\cdot p_2^*o_\mathcal{X})]\,\,\,
\\
\nonumber
-[ Z_2\circ (p_2^*o_\mathcal{X}\cdot p_3^*o_\mathcal{X})]
+[ Z_4\circ (p_1^*o_\mathcal{X}\cdot p_3^*o_\mathcal{X})]
+[ Z_4\circ (p_2^*o_\mathcal{X}\cdot p_3^*o_\mathcal{X})],\,\,\,
\end{eqnarray}
which rewrites as
\begin{eqnarray}\label{equanou4aout2011matin}
[\Delta]=[ Z_0\circ (p_1^*o_\mathcal{X}\cdot p_2^*o_\mathcal{X})]+[ Z_4\circ \Delta]\,\,\,\,\,\,\\
\nonumber
-[ Z_4\circ (p_{13}^*\Delta_\mathcal{X}\cdot p_2^*o_\mathcal{X})]
-[ Z_4\circ (p_1^*o_\mathcal{X}\cdot p_{23}^*\Delta_\mathcal{X})]+[Z_4\circ (p_1^*o_\mathcal{X}\cdot p_3^*o_\mathcal{X})]
\,\,\,\,\,\,
\\
\nonumber
+[Z_4\circ (p_2^*o_\mathcal{X}\cdot p_3^*o_\mathcal{X})]
+[Z_4\circ (p_1^*o_\mathcal{X}\cdot p_2^*o_\mathcal{X})]
+[ Z_2\circ (p_1^*o_\mathcal{X}\cdot p_{23}^*\Delta_\mathcal{X})]\\
\nonumber
\,\,\,\,\,\,
-[ Z_2\circ (p_1^*o_\mathcal{X}\cdot p_3^*o_\mathcal{X})]
+[ Z_2\circ (p_{13}^*\Delta_\mathcal{X}\cdot p_2^*o_\mathcal{X})]
-2[ Z_2\circ (p_1^*o_\mathcal{X}\cdot p_2^*o_\mathcal{X})].
\,\,\,\,\,\,
\end{eqnarray}
To conclude, we use the following lemma:
\begin{lemm}  \label{lemmesameacabit}Up to passing to a dense Zariski open set
of $B$, we have the following equalities in $CH^4(\mathcal{X}\times_B\mathcal{X}\times_B\mathcal{X})_\mathbb{Q}$:
\begin{eqnarray}\label{eq4aoutmarre0}Z_0\circ (p_1^*o_\mathcal{X}\cdot p_2^*o_\mathcal{X})=p_1^*o_\mathcal{X}\cdot p_2^*o_\mathcal{X},
\end{eqnarray}
\begin{eqnarray}\label{eq4aoutmarre1} Z_4\circ \Delta=p_{12}^*\Delta_\mathcal{X}\cdot p_3^*o_\mathcal{X},
\end{eqnarray}
\begin{eqnarray}\label{eq4aoutmarre2} Z_4\circ (p_{13}^*\Delta_\mathcal{X}\cdot p_2^*o_\mathcal{X})=p_2^*o_\mathcal{X}\cdot p_3^*o_\mathcal{X},
\end{eqnarray}
\begin{eqnarray}\label{eq4aoutmarre3} Z_4\circ (p_2^*o_\mathcal{X}\cdot p_3^*o_\mathcal{X})=p_2^*o_\mathcal{X}\cdot p_3^*o_\mathcal{X},
\end{eqnarray}
\begin{eqnarray}\label{eq4aoutmarrebisbis} Z_4\circ (p_1^*o_\mathcal{X}\cdot p_{23}^*\Delta_\mathcal{X})=p_1^*o_\mathcal{X}\cdot p_3^*o_\mathcal{X},
\end{eqnarray}
\begin{eqnarray}\label{eq4aoutmarrebister}Z_4\circ (p_1^*o_\mathcal{X}\cdot p_2^*o_\mathcal{X})=0,
\end{eqnarray}
\begin{eqnarray}\label{eq4aoutmarrebisterter}Z_4\circ (p_1^*o_\mathcal{X}\cdot p_3^*o_\mathcal{X})=p_1^*o_\mathcal{X}\cdot p_3^*o_\mathcal{X},
\end{eqnarray}
\begin{eqnarray}\label{eq4aoutmarre7} Z_2\circ (p_1^*o_\mathcal{X}\cdot p_3^*o_\mathcal{X})=0,
\end{eqnarray}
\begin{eqnarray}\label{eq4aoutmarre9}  Z_2\circ (p_1^*o_\mathcal{X}\cdot p_{23}^*\Delta_\mathcal{X})=p_1^*o_\mathcal{X}\cdot p_{23}^*\Delta_\mathcal{X}-p_1^*o_\mathcal{X}\cdot p_2^*o_\mathcal{X}-p_1^*o_\mathcal{X}\cdot p_3^*o_\mathcal{X},
\end{eqnarray}
\begin{eqnarray}\label{eq4aoutmarre9bis}
Z_2\circ (p_{13}^*\Delta_\mathcal{X}\cdot p_2^*o_\mathcal{X})=p_{13}^*\Delta_\mathcal{X}\cdot p_2^*o_\mathcal{X}-p_1^*o_\mathcal{X}\cdot p_2^*o_\mathcal{X}-
p_2^*o_\mathcal{X}\cdot p_3^*o_\mathcal{X},
\end{eqnarray}
\begin{eqnarray}\label{eq4aoutmarre8} Z_2\circ (p_1^*o_\mathcal{X}\cdot p_2^*o_\mathcal{X})=0.
\end{eqnarray}

\end{lemm}
{\bf Proof.} The proof of (\ref{eq4aoutmarre2}) is explicit, recalling that
$Z_4=\{(x,o_{\mathcal{X}_b}),\,x\in \mathcal{X}_b,\,b\in B\}$, and that
$p_{13}^*\Delta_\mathcal{X}\cdot p_2^*o_\mathcal{X}=\{(y,o_{\mathcal{X}_b},y),\,y\in \mathcal{X}_b,\,b\in B\}$.
We then find that
$Z_4\circ (p_{13}^*\Delta_\mathcal{X}\cdot p_2^*o_\mathcal{X})$ is the cycle
$$ p_{124}(p_{13}^*\Delta_\mathcal{X}\cdot p_2^*o_\mathcal{X}\cdot p_{34}^*(Z_4))=
p_{124}(\{(y,o_{\mathcal{X}_b},y,o_{\mathcal{X}_b}),\,y\in \mathcal{X}_b,\,b\in B\})$$
$$=\{(y,o_{\mathcal{X}_b},o_{\mathcal{X}_b}),\,y\in \mathcal{X}_b,\,b\in B\},$$
which proves (\ref{eq4aoutmarre2}).
(\ref{eq4aoutmarre3}) is the same formula as (\ref{eq4aoutmarre0}) with the indices $1$ and $3$
exchanged. The proofs of
(\ref{eq4aoutmarre0}) to (\ref{eq4aoutmarrebisterter}) work similarly.

For the other proofs, we recall that
$$Z_2=\Delta_\mathcal{X}-Z_0-Z_4\subset \mathcal{X}\times_B\mathcal{X}.$$
Thus we get, as $\Delta_\mathcal{X}$ acts as the identity:
$$ Z_2\circ (p_1^*o_\mathcal{X}\cdot p_{23}^*\Delta_\mathcal{X})=p_1^*o_\mathcal{X}\cdot p_{23}^*\Delta_\mathcal{X}-Z_0\circ (p_1^*o_\mathcal{X}\cdot p_{23}^*\Delta_\mathcal{X})-Z_4\circ(p_1^*o_\mathcal{X}\cdot p_{23}^*\Delta_\mathcal{X}).$$
We then compute the terms
$Z_0\circ (p_1^*o_\mathcal{X}\cdot p_{23}^*\Delta_\mathcal{X}),\,Z_4\circ(p_1^*o_\mathcal{X}\cdot p_{23}^*\Delta_\mathcal{X})$ explicitly as before, which gives (\ref{eq4aoutmarre9}).

The other  proofs are similar.

\cqfd

Using the cohomological version of Lemma \ref{lemmesameacabit}, (\ref{equanou4aout2011matin}) becomes:
\begin{eqnarray}\label{catastrophique}
[\Delta]=[ p_1^*o_\mathcal{X}\cdot p_2^*o_\mathcal{X}]+[p_{12}^*\Delta_\mathcal{X}\cdot p_3^*o_\mathcal{X}]\,\,\,\,\,\,\\
\nonumber
-[p_2^*o_\mathcal{X}\cdot p_3^*o_\mathcal{X})]-[p_1^*o_\mathcal{X}\cdot p_3^*o_\mathcal{X})]
+[p_1^*o_\mathcal{X}\cdot p_3^*o_\mathcal{X}]\,\,\,\,\,\,
\\
\nonumber+[p_2^*o_\mathcal{X}\cdot p_3^*o_\mathcal{X}]
+[ p_1^*o_\mathcal{X}\cdot p_{23}^*\Delta_\mathcal{X}-p_1^*o_\mathcal{X}\cdot p_2^*o_\mathcal{X}-p_1^*o_\mathcal{X}\cdot p_3^*o_\mathcal{X}]
\\
\nonumber
\,\,\,\,\,\,
+[ p_{13}^*\Delta_\mathcal{X}\cdot p_2^*o_\mathcal{X}-p_1^*o_\mathcal{X}\cdot p_2^*o_\mathcal{X}-
p_2^*o_\mathcal{X}\cdot p_3^*o_\mathcal{X}].
\,\,\,\,\,\,
\end{eqnarray}
This last equality is now satisfied by assumption (compare with (\ref{equasmalldiagrel29ju})) and this concludes the proof
of formula (\ref{equanou3aout2011}). Thus (i) is proved.

\vspace{0.5cm}

(ii)  We just have to prove that
\begin{eqnarray}\label{annu29ju1438}
P_0\otimes P_0 ([\Delta_\mathcal{X}])=P_4\otimes P_4 ([\Delta_\mathcal{X}])=0,\\
\nonumber
 P_0\otimes P_2([\Delta_\mathcal{X}])=P_4\otimes P_2([\Delta_\mathcal{X}])=0\,\,
{\rm in}\,\, H^4(\mathcal{X}\times_B\mathcal{X},\mathbb{Q}).
\end{eqnarray}
Indeed, the relative K\"unneth decomposition gives
$$ R(\pi,\pi)_*\mathbb{Q}=R\pi_*\mathbb{Q}\otimes R\pi_*\mathbb{Q}$$
and
the decomposition
(\ref{ladecomp29ju2011}) induces a
decomposition of the above tensor product on the right:
\begin{eqnarray}\label{cohindecomp29ju}
R\pi_*\mathbb{Q}\otimes R\pi_*\mathbb{Q}=\oplus _{k,l}R^k\pi_*\mathbb{Q}\otimes R^l\pi_*\mathbb{Q}[-k-l],
\end{eqnarray}
where the decomposition is induced by the various tensor products
of $P_0,P_2,P_4$.
Taking cohomology
in (\ref{cohindecomp29ju}) gives
$$H^4(\mathcal{X}\times_B\mathcal{X},\mathbb{Q})=\oplus_{s+k+l=4}H^s(B,R^k\pi_*\mathbb{Q}\otimes R^l\pi_*\mathbb{Q}).$$
The term
$H^0(R^4(\pi,\pi)_*\mathbb{Q})$ is then exactly the term in the above decomposition
of $H^4(\mathcal{X}\times_B\mathcal{X},\mathbb{Q})$ which is annihilated by the four projectors
$P_0\otimes P_0 $, $P_0\otimes P_2$, $P_4\otimes P_2$, $P_4\otimes P_4$ and those obtained by changing the order of factors.

The proof of (\ref{annu29ju1438}) is elementary. Indeed, consider for example
the term $P_0\otimes P_0$, which is given by the cohomology class of the cycle
$$Z:=pr_1^*o_\mathcal{X}\cdot pr_2^*o_\mathcal{X}\subset \mathcal{X}\times_B\mathcal{X}\times_B\mathcal{X}\times_B\mathcal{X},$$
which we see as a relative self-correspondence of
$\mathcal{X}\times_B\mathcal{X}$
We have
$$Z_*(\Delta_\mathcal{X})= p_{34*}(p_{12}^*\Delta_\mathcal{X}\cdot Z).$$
But the cycle on the right is trivially rationally equivalent to $0$ on fibers
$\mathcal{X}_t\times \mathcal{X}_t$. It thus follows from the general principle \ref{genprinciple} that for some dense Zariski open set $B^0$ of $B$,
$$[Z]_*([\Delta_\mathcal{X}])=0\, \,{\rm in}\,\,H^4(\mathcal{X}^0\times_{B^0}\mathcal{X}^0,\mathbb{Q}).$$
The other vanishing statements  are proved similarly.
\cqfd

\section{Calabi-Yau hypersurfaces \label{section3}}
In the case of smooth Calabi-Yau hypersurfaces $X$ in projective space
$\mathbb{P}^n$, that is hypersurfaces of degree $n+1$ in $\mathbb{P}^n$, we have
the following result which partially generalizes Theorem
 \ref{decompdiagonalpetiteintro} and
 provides some information on the Chow ring of $X$. Denote by $o\in CH_0(X)_\mathbb{Q}$ the class of the $0$-cycle
 $\frac{h^{n-1}}{n+1}$, where $h:=c_1(\mathcal{O}_X(1))\in CH^1(X)$. We denote again by
 $\Delta$ the small diagonal of $X$ in $X^3$.
 \begin{theo} \label{theonewpourarticle} The following relation  holds in $CH^{2n-2}(X\times X\times X)_\mathbb{Q}$:
 \begin{eqnarray}\label{equadeltapetiteCY}\Delta=\Delta_{12}\cdot o_3+(perm.)+Z+\Gamma',
\end{eqnarray}
where $Z$ is the restriction to $X\times X\times X$
of a cycle on $\mathbb{P}^n\times\mathbb{P}^n\times\mathbb{P}^n$, and $\Gamma'$ is a multiple
of the following effective cycle of dimension $n-1$:
\begin{eqnarray}
\label{nouveaunom10aout}\Gamma:=\cup_{t\in F(X)}\mathbb{P}^1_t\times\mathbb{P}^1_t\times\mathbb{P}^1_t.
\end{eqnarray}
 \end{theo}

 Here $F(X)$ is the variety of lines contained in $X$. It is of dimension $n-4$ for general $X$. For $t\in F(X)$ we denote
 $\mathbb{P}^1_t\subset X\subset \mathbb{P}^n$ the corresponding line.

 \vspace{0.5cm}

 {\bf Proof of Theorem \ref{theonewpourarticle}.} Observe first of all that it suffices to prove
 the following equality of $n-1$-cycles on  $X^3_0:=X^3\setminus\Delta$:
  \begin{eqnarray}\label{equadeltapetiteCY1proof}\Gamma_{\mid X^3_0}=2(n+1)(\Delta_{12\mid X^3_0}\cdot o_3+(perm.))+Z\,\,{\rm in}\,\,CH^{2n-2}(X^3_0)_\mathbb{Q}
\end{eqnarray}
where $Z$ is the restriction to $X^3_0$
of a cycle on $(\mathbb{P}^n)^3$. Indeed, by the localization exact sequence
(cf. \cite[Lemma 9.12]{voisinbook}), (\ref{equadeltapetiteCY1proof}) implies an equality, for an adequate multiple $\Gamma'$ of $\Gamma$:
\begin{eqnarray}\label{equadeltapetiteCYproof2}N\Delta=\Delta_{12}\cdot o_3+(perm.)+Z+\Gamma'\,\,{\rm in}\,\,CH^{2n-2}(X\times X\times X)_\mathbb{Q},
\end{eqnarray}
for some rational number $N$.
Projecting to $X^2$ and taking cohomology classes, we easily conclude then that $N=1$.
(We use here the fact that $X$ has some transcendental cohomology, so that the
cohomology  class of
the diagonal of $X$ does not vanish on products $U\times U$, where $U\subset X$
is Zariski open.)

In order to prove (\ref{equadeltapetiteCY1proof}), we do the following:
First of all we compute the class in $CH^{n-1}(X^3_0)$ of
the $2n-2$-dimensional subvariety
$$X^3_{0,col, sch}\subset X^3_0$$
parameterizing $3$-uples of collinear points satisfying the following property:

{\it Let $\mathbb{P}^1_{x_1x_2x_3}=<x_1,x_2,x_3>$ be the line generated by the
$x_i$'s. Then the  subscheme $x_1+x_2+x_3$ of $\mathbb{P}^1_{x_1x_2x_3}\subset \mathbb{P}^n$ is contained in $X$.}

We will denote $X^3_{0,col}\subset X^3_0$
the $2n-2$-dimensional subvariety parameterizing $3$-uples of collinear points.
 Obviously
 $X^3_{0,col, sch}\subset X^3_{0,col}$.
 We will see that the first one is in fact an irreducible component of the second one.

Next we observe that  there
is a natural morphism
 $\phi: X^3_{0,col}\rightarrow G(2,n+1)$
 to the Grassmannian  of lines in $\mathbb{P}^n$,
which to $(x_1,x_2,x_3)$ associates the line $\mathbb{P}^1_{x_1x_2x_3}$. This morphism is well-defined
on $X^3_{0,col}$ because at least two of the points $x_i$ are distinct, so that
this line is well-determined.
The morphism $\phi$ corresponds to a tautological rank $2$ vector bundle
$\mathcal{E}$ on $X^3_{0,col}$, with fiber $H^0(\mathcal{O}_{\mathbb{P}^1_{x_1x_2x_3}}(1))$ over the point
$(x_1,x_2,x_3)$.
We then observe that $\Gamma\subset X^3_{0,col, sch}$ is defined by the condition
that the line $\mathbb{P}^1_{x_1x_2x_3}$ be contained in $X$. In other words, the
equation $f$ defining $X$ has to vanish on this line. This condition can be seen globally as the vanishing of the section
$\sigma$ of
the vector bundle $S^{n+1}\mathcal{E}$ defined by
$$\sigma((x_1,x_2,x_3))=f_{\mid\mathbb{P}^1_{x_1x_2x_3}},$$

This section $\sigma$ is not transverse, (in fact the rank of $S^{n+1}\mathcal{E}$ is
$n+2$, while the codimension of $\Gamma$ is $n-1$), but the reason for this is very simple:
indeed, at a point $(x_1,x_2,x_3)$ of $X^3_{0,col,sch}$, the equation
$f$ vanishes by definition on the degree $3$ cycle $x_1+x_2+x_3$ of
$\mathbb{P}^1_{x_1x_2x_3}$.
Another way to express this is to say that
$\sigma$ is in fact a section of the rank $n-1$ bundle
\begin{eqnarray}\label{vbundle25juillet}\mathcal{F}\subset S^{n+1}\mathcal{E}
\end{eqnarray}
where $\mathcal{F}_{(x_1,x_2,x_3)}$ consists of degree $n+1$ polynomials
vanishing on the { subscheme} $x_1+x_2+x_3$ of $\mathbb{P}^1_{x_1x_2x_3}$.

The section $\sigma$
of $\mathcal{F}$ is transverse and thus we conclude that we have the following equality
\begin{eqnarray}\label{equadeltapetiteCYproof3}\Gamma_{\mid X^3_0}=j_*(c_{n-1}(\mathcal{F}))\,\,{\rm in}\,\,CH^{2n-2}(X^3_0)_\mathbb{Q},
\end{eqnarray}
where $j$ is the inclusion of $X^3_{0,col,sch}$ in $X^3_0$.

We now observe that the vector bundles  $\mathcal{E}$ and $\mathcal{F}$
come from  vector bundles
on the variety $(\mathbb{P}^n)^3_{0,col}$ parameterizing $3$-uples of
collinear points in $\mathbb{P}^n$, at least two of them being distinct.

The variety $(\mathbb{P}^n)^3_{0,col}$ is smooth irreducible of dimension
 $2n+1$ (hence of codimension $n-1$ in $(\mathbb{P}^n)^3$), being Zariski open
in a $\mathbb{P}^1\times \mathbb{P}^1\times \mathbb{P}^1$-bundle over the Grassmannian
$G(2,n+1)$.
We have now the following:
\begin{lemm}
The intersection
$(\mathbb{P}^n)^3_{0,col}\cap X^3_0$
is reduced, of pure dimension $2n-2$. It decomposes as
\begin{eqnarray}\label{decompinters}
(\mathbb{P}^n)^3_{0,col}\cap X^3_0=X^3_{0,col,sch}\cup \Delta_{0, 12}\cup \Delta_{0, 13}\cup \Delta_{0, 23},
\end{eqnarray}
where $\Delta_{0,ij}\subset X^3_0$ is defined as  ${\Delta_{ij}}\cap X_0^3$
with $\Delta_{ij}$ the big diagonal $\{x_i=x_j\}$.
\end{lemm}
{\bf Proof.} The set theoretic equality in (\ref{decompinters}) is obvious.
The fact that each component on the right has dimension $2n-2$ and thus is a component
of the right dimension of this intersection is also obvious. The only
point to check is thus the fact that these intersections are transverse at the generic point of each
component in the right hand side.
The generic point of the irreducible variety
$X^3_{0,col,sch}$ parameterizes  a triple of distinct collinear points which are on a
line $D$ not tangent to $X$. At such a triple,
the intersection $(\mathbb{P}^n)^3_{0,col}\cap X^3_0$ is smooth of dimension
$2n-2$ because $(\mathbb{P}^n)^3_{0,col}$ is Zariski open in the triple self-product
$P\times_{G(2,n+1)} P\times_{G(2,n+1)}P$ of the tautological $\mathbb{P}^1$-bundle $P$ over the Grassmannian
$G(2,n+1)$, and the intersection with $X^3_0$ is defined by the three equations
$$p\circ pr_1^*f,\,p\circ pr_2^*f,\,p\circ pr_3^*f,$$
where the $pr_i$'s are the projections $P^{3/G(2,n+1)}\rightarrow P$
 and $p:P\rightarrow \mathbb{P}^n$ is the natural map.
These three equations are independent since they are independent
after restriction to  $D\times D\times D \subset P\times_{G(2,n+1)} P\times_{G(2,n+1)}P$
at the point
$(x_1,x_2,x_3)$ because $D$ is not tangent to $X$.

Similarly, the generic point of the irreducible variety
$\Delta_{0,1,2,j}\subset X^3_{0,col}$ parameterizes a triple
$(x,x,y)$ with the property that $x\not=y$ and the line $\mathbb{P}^1_{xy}:=<x,y>$ is
not tangent to $X$. Again, the intersection
$(\mathbb{P}^n)^3_{0,col}\cap X^3_0$ is smooth of dimension
$2n-2$ near $(x,x,y)$ because the restrictions to
$\mathbb{P}^1_{xy}\times \mathbb{P}^1_{xy}\times \mathbb{P}^1_{xy}\subset P\times_{G(2,n+1)} P\times_{G(2,n+1)}P$ of
the equations  $$p\circ pr_1^*f,\,p\circ pr_2^*f,\,p\circ pr_3^*f,$$
defining $X^3$ are independent.

\cqfd

Combining (\ref{decompinters}), (\ref{equadeltapetiteCYproof3}) and the fact that
the vector bundle $\mathcal{F}$ already exists on $(\mathbb{P}^n)^3_{0,col}$,
we find that
$$\Gamma_{\mid X^3_0}=J_*(c_{n-1}(\mathcal{F}_{\mid (\mathbb{P}^n)^3_{0,col}\cap X^3_0}))-
\sum_{i\not=j}J_{0,ij*}c_{n-1}(\mathcal{F}_{\mid \Delta_{0,ij}})\,\,{\rm in}\,\,CH^{2n-2}(X^3_0)_\mathbb{Q},
$$
where $J:(\mathbb{P}^n)^3_{0,col}\cap X^3_0\hookrightarrow  X^3_0$
is the inclusion and similarly for $J_{0,ij}: \Delta_{0,ij}\hookrightarrow  X^3_0$.
This provides us with the formula:
\begin{eqnarray}\label{equadeltapetiteCYproof4}\Gamma_{\mid X^3_0}=(K_*c_{n-1}(\mathcal{F}))_{\mid X^3_0} -
\sum_{i\not=j}J_{0,ij*}c_{n-1}(\mathcal{F}_{\mid \Delta_{0,ij}})\,\,{\rm in}\,\,CH^{2n-2}(X^3_0)_\mathbb{Q},
\end{eqnarray}
where $K: (\mathbb{P}^n)^3_{0,col}\hookrightarrow (\mathbb{P}^n)^3_{0}$ is the inclusion map.

The first term comes from $CH((\mathbb{P}^n)^3_{0})$, so to conclude we only have to compute
the terms $J_{0,ij*}c_{n-1}(\mathcal{F}_{\mid \Delta_{0,ij}})$.
This is however very easy, because the vector bundles $\mathcal{E}$ and $\mathcal{F}$ are very simple
on $\Delta_{0,ij}$: Assume for simplicity $i=1,j=2$. Points of $\Delta_{0,12}$ are points
$(x,x,y),\,x\not= y\in X$.
The line $\phi((x,x,y))$ is the line $<x,y>,\,x\not=y$, and it follows that
\begin{eqnarray}
\label{finequa25juillet}\mathcal{E}_{\mid \Delta_{0,12}}=pr_2^*\mathcal{O}_X(1)\oplus pr_3^*\mathcal{O}_X(1).
\end{eqnarray}
The projective bundle $\mathbb{P}(\mathcal{E}_{\mid \Delta_{0,12}})$ has two sections on
$\Delta_{0,12}$ which give two divisors
$$D_2\in|\mathcal{O}_{\mathbb{P}(\mathcal{E})}(1)\otimes pr_3^*\mathcal{O}_X(-1)|,\,
D_3\in |\mathcal{O}_{\mathbb{P}(\mathcal{E})}(1)\otimes pr_2^*\mathcal{O}_X(-1)|.$$
The length $3$ subscheme $2D_2+D_3\subset \mathbb{P}(\mathcal{E}_{\mid \Delta_{0,1,2}})$
with fiber $2x+y$ over the point $(x,x,y)$  is thus
the zero set of
a section $\alpha$ of the line bundle
$\mathcal{O}_{\mathbb{P}(\mathcal{E})}(3)\otimes pr_3^*\mathcal{O}_X(-2)\otimes pr_2^*\mathcal{O}_X(-1)$.
We thus conclude that the vector bundle $\mathcal{F}_{\mid \Delta_{0,12}}$ is
isomorphic to
$$pr_3^*\mathcal{O}_X(2)\otimes pr_2^*\mathcal{O}_X(1)\otimes S^{n-2}\mathcal{E}_{\mid \Delta_{0,12}}.$$
Combining with (\ref{finequa25juillet}), we conclude that
$c_{n-1}(\mathcal{F}_{\mid \Delta_{0,12}})$ can be expressed as a polynomial of
degree $n-1$ in
$h_2=c_1(pr_2^*\mathcal{O}_X(1))$ and $h_3=c_1(pr_3^*\mathcal{O}_X(1)))$ on
$\Delta_{0,12}$.
The proof of (\ref{equadeltapetiteCY1proof}) is completed by the following lemma:
\begin{lemm} Let $\Delta_X\subset X\times X$ be the diagonal. Then
the codimension $n$ cycles
$$ pr_1^*c_1(\mathcal{O}_X(1))\cdot \Delta_X,\,\,pr_2^*c_1(\mathcal{O}_X(1))\cdot \Delta_X$$
of $X\times X$
are restrictions to $X\times X$ of cycles $Z\in CH^n(\mathbb{P}^n\times \mathbb{P}^n)_\mathbb{Q}$.
\end{lemm}
{\bf Proof.}  Indeed, let $j_X:X\hookrightarrow \mathbb{P}^n$ be the inclusion
of $X$ in $\mathbb{P}^n$, and $j_{X,1},\,j_{X,2}$ the corresponding inclusions
of $X\times X$ in $\mathbb{P}^n\times X$, resp. $X\times \mathbb{P}^n$.
Then as $X$ is a degree $n+1$ hypersurface, the composition $j_{X,1}^*\circ  j_{X,1*}: CH^*(X\times X)\rightarrow CH^{*+1}(X\times X)$
is equal to the morphism given by intersection with the class $(n+1)pr_1^*c_1(\mathcal{O}_X(1))$,
and similarly for the second inclusion.
On the other hand, $j_{X,1*}(\Delta_X)\subset \mathbb{P}^n\times X$ is obviously the
(transpose of the)  graph of the
inclusion of $X$ in $\mathbb{P}^n$, hence its class is the restriction to
$\mathbb{P}^n\times X$ of the diagonal of $\mathbb{P}^n\times \mathbb{P}^n$. This implies that
$$(n+1)pr_1^*c_1(\mathcal{O}_X(1))\cdot \Delta_X=j_{X,1}^*((\Delta_{\mathbb{P}^n\times \mathbb{P}^n})_{\mid \mathbb{P}^n\times X}),$$
which proves the result for $ pr_1^*c_1(\mathcal{O}_X(1))\cdot \Delta_X$.
We argue similarly for the second cycle.
\cqfd
It follows from this lemma that
a monomial of
degree $n-1$ in
$h_2=c_1(pr_2^*\mathcal{O}_X(1))$ and $h_3=c_1(pr_3^*\mathcal{O}_X(1)))$ on
$\Delta_{0,12}$, seen as a cycle in $X^3_0$, will be the restriction
to $X^3_0$ of a cycle with $\mathbb{Q}$-coefficients on $(\mathbb{P}^n)^3$, unless it is
proportional to $h_3^{n-1}$.
Recalling that $c_1(\mathcal{O}_X(1))^{n-1}=(n+1)o\in CH_0(X)$, we finally
proved that modulo
restrictions of cycles on
$(\mathbb{P}^n)^3$, the term
$J_{0,12*}c_{n-1}(\mathcal{F}_{\mid \Delta_{0,12}})$ is a multiple of
$(\Delta_{12}\cdot o_3)_{\mid X^3_0}$ in $CH^{2n-2}(X^3_0)_\mathbb{Q})$.
 The precise coefficient is in fact given by the argument above. Indeed,
 we just saw that modulo restrictions of cycles coming from
$\mathbb{P}^n\times \mathbb{P}^n\times \mathbb{P}^n$, the term
$J_{0,12*}c_{n-1}(\mathcal{F}_{\mid \Delta_{0,12}})$ is equal to
\begin{eqnarray}\label{derproofCY26}\mu \Delta_{12}\cdot pr_3^*(c_1(\mathcal{O}_X(1))^{n-1})=\mu(n+1)(\Delta_{12}\times o_3)_{\mid X^3_0},
\end{eqnarray}
with $c_1(\mathcal{O}_X(1))^{n-1}=(n+1)o$ in $CH_0(X)$, and where the coefficient $\mu$
is the coefficient of $h_3^{n-1}$ in the polynomial in $h_2,\,h_3$ computing
$c_{n-1}(\mathcal{F}_{\mid \Delta_{0,12}})$.

We use now the isomorphism
$$\mathcal{F}_{\mid \Delta_{0,12}}\cong pr_3^*\mathcal{O}_X(2)\otimes pr_2^*\mathcal{O}_X(1)\otimes S^{n-2}\mathcal{E}_{\mid \Delta_{0,12}},$$
where
$\mathcal{E}_{\mid \Delta_{0,12}}\cong pr_2^*\mathcal{O}_X(1)\oplus pr_3^*\mathcal{O}_X(1)$
according to (\ref{finequa25juillet}).
Hence we conclude that the coefficient $\mu$ is equal to
$2$, and this concludes the proof
 of (\ref{equadeltapetiteCY1proof}), using (\ref{derproofCY26}) and (\ref{equadeltapetiteCYproof4}).

 \cqfd
 We have the following consequence of Theorem \ref{theonewpourarticle}, which is a generalization of
 Theorem \ref{BV} to Calabi-Yau hypersurfaces.
 \begin{theo} \label{theogenBV} Let $Z_i,\,Z'_i$ be cycles of codimension $>0$
 on $X$ such that ${\rm codim}\,Z_i+{\rm codim}\,Z'_i= n-1$. Then if we have a cohomological relation
 $$\sum_i n_i[Z_i]\cup[Z'_i]=0\,\,{\rm in}\,\,H^{2n-2}(X,\mathbb{Q})$$
 this relation already holds at the level of Chow groups:
 $$\sum_i n_iZ_i\cdot Z'_i=0\,\,{\rm in}\,\,CH_0(X)_\mathbb{Q}.$$
 \end{theo}
 {\bf Proof.} Indeed, let us view formula (\ref{equadeltapetiteCY})
 as an equality of correspondences between $X\times X$ and $X$.
 The left hand side applied to $\sum_i n_iZ_i\times Z'_i$
 is the desired cycle:
 $\Delta_*(\sum_i n_iZ_i\times Z'_i)=\sum_i n_iZ_i\cdot Z'_i\,\,{\rm in}\,\,CH_0(X)_\mathbb{Q}$.
 The right hand side is a sum of three terms:
 \begin{eqnarray}\label{cycle12aout}
 (\Delta_{12}\cdot o_3+(perm.))_*(\sum_i n_iZ_i\times Z'_i)+Z_*(\sum_i n_iZ_i\times Z'_i)+\Gamma'_*(\sum_i n_iZ_i\times Z'_i).
 \end{eqnarray}
 For the first term, we observe that  $(\Delta_{12}\cdot o_3)_*(\sum_i n_iZ_i\times Z'_i)=({\rm deg}\,\sum_i n_iZ_i\cdot Z'_i)\,o_3$
 vanishes in $CH_0(X)_\mathbb{Q}$, and that the two other terms
 $(\Delta_{13}\cdot o_2)_*(\sum_i n_iZ_i\times Z'_i)$ and $(\Delta_{23}\cdot o_1)_*(\sum_i n_iZ_i\times Z'_i)$
 vanish by the assumption that ${\rm codim}\,Z_i>0$ for all $i$.

 For the second term,  we observe that as $Z$ is the restriction
 of a cycle $Z'\in CH^{2n-2}(\mathbb{P}^n\times\mathbb{P}^n\times\mathbb{P}^n)_\mathbb{Q}$,
 $Z_*(\sum_i n_iZ_i\times Z'_i)$ is equal to
 $$j^*({Z'}^*((j,j)_*(\sum_i n_iZ_i\times Z'_i)))\in CH^{n-1}(X)_\mathbb{Q}.$$
 Hence it belongs to ${\rm Im}\,j^*$, and is proportional to $o$.

 Consider finally the term $\Gamma'_*(\sum_i n_iZ_i\times Z'_i)$, which is a multiple of
 $\Gamma'_*(\sum_i n_iZ_i\times Z'_i)$: Let
 $\Gamma_0\subset X$ be the locus swept-out by lines. We observe that
 for any line $D\cong \mathbb{P}^1\subset X$, any point on $D$ is rationally equivalent to
 the zero cycle $h\cdot D$ which is in fact proportional to $o$, since
 $$(n+1) h\cdot D=j^*\circ j_*(D)\,\,{\rm in}\,\, CH_0(X)$$
 and $j_*(D)=c_1(\mathcal{O}_{\mathbb{P}^n}(1))^{n-1}\,\,{\rm in}\,\, CH^{n-1}(\mathbb{P}^n)$.
 Hence all points of $\Gamma_0$ are rationally equivalent  to $o$ in $X$, and thus
 $\Gamma'_*(\sum_i n_iZ_i\times Z'_i)$ is also proportional to $o$.

 It follows from the above analysis that  the $0$-cycle (\ref{cycle12aout}) is
  a multiple of $o$ in $CH_0(X)_\mathbb{Q}$. As it is of degree $0$, it is in fact rationally equivalent to $0$.

 \cqfd
 We leave as a conjecture the following :
 \begin{conj}\label{final17janv} For any smooth $n-1$-dimensional Calabi-Yau hypersurface
for which the variety of lines  $F(X)$ has dimension $n-4$, the $n-1$-cycle $\Gamma\in CH_{n-1}( X\times X\times X)_\mathbb{Q}$
 of (\ref{nouveaunom10aout}) is the restriction to $X\times X\times X$ of a $n+2$-cycle on $\mathbb{P}^n\times\mathbb{P}^n\times\mathbb{P}^n$.
\end{conj}


\begin{thebibliography}{99}
\bibitem{beauville} A. Beauville. On the splitting of the Bloch-Beilinson filtration, in {\it Algebraic cycles and motives} (vol. 2), London Math. Soc. Lecture Notes 344, 38-53, Cambridge University Press (2007).
\bibitem{beauvoi} A. Beauville, C. Voisin. On the Chow ring of a $K3$ surface,  J. Algebraic Geom.  13  (2004),  no. 3, 417-426.
\bibitem{blochsrinivas} S. Bloch and V. Srinivas. Remarks on
 correspondences and algebraic cycles,
 Amer. J. of Math. 105 (1983) 1235-1253.
 \bibitem{deligne} P. Deligne.  Th\'{e}or\`{e}me de Lefschetz et crit\`{e}res de d\'{e}g\'{e}n\'{e}rescence de
suites spectrales,  Inst. Hautes Etudes Sci. Publ. Math. No. 35,
 259-278, (1968).
 \bibitem{delignemotivedec} P. Deligne. D\'{e}compositions dans la cat\'{e}gorie
d\'{e}riv\'{e}e, Motives (Seattle, WA, 1991), 115-128, Proc. Sympos. Pure Math., 55, Part 1.
\bibitem{deningermurre} C. Deninger, J. P. Murre. Motivic decomposition of abelian schemes and the Fourier transform. J. Reine Angew. Math. 422 (1991), 201-219.
    \bibitem{fulton} W. Fulton. {\it Intersection Theory}, Ergebnisse der Math. und ihrer Grenzgebiete 3 Folge, Band 2, Springer (1984).
    \bibitem{godement} R. Godement. {\it Topologie alg\'{e}brique et th\'{e}orie des faisceaux},  Actualit\'{e}s scientifiques et industrielles 1252, Hermann, Paris (1964).
 \bibitem{huy} D. Huybrechts. A note on the Bloch-Beilinson conjecture for $K3$ surfaces and spherical objects, Pure and Applied
Mathematics Quarterly, Special volume in Memory of Eckart Viehweg,
Volume 7, Number 4,	p.1395-1406
 \bibitem{mumford} D. Mumford. Rational equivalence of zero-cycles on surfaces,
J. Math. Kyoto Univ. 9 (1968), 195-204.
\bibitem{murre} J. P. Murre. On the motive of an algebraic surface.  J. Reine Angew. Math.  409  (1990), 190-204.
\bibitem{voisinKcorresp} C. Voisin. Intrinsic pseudo-volume forms and $K$-correspondences, in  {\it The Fano Conference},  761-792, Univ. Torino, Turin, (2004).
    \bibitem{voisin} C. Voisin.  On the Chow ring of certain algebraic hyper-K\"{a}hler manifolds, Pure and Applied Mathematics Quarterly, Vol. 4, N0 3 613-649 (2008).
   \bibitem{voisinbook} C. Voisin. {\it Hodge theory and complex algebraic geometry} II,  Cambridge Studies in Advanced Mathematics, 77. Cambridge University Press, Cambridge, (2003).
       \bibitem{voisinweyllectures} C. Voisin. {\it Chow rings, decomposition of the diagonal and the topology of families}, Hermann Weyl lectures 2011, to appear in the Annals of Math. Studies.
\end{thebibliography}
\end{document}